\documentclass[12pt,reqno]{amsart}
\usepackage{amsmath,amssymb,enumerate, color}
\usepackage[dvipdfmx]{graphicx}
\usepackage[all]{xy}
\usepackage[varg]{txfonts}
\usepackage{bm}

\newtheorem{tm}{Theorem}[section]
\newtheorem{lm}[tm]{Lemma}

\newtheorem{re}[tm]{Remark}

\newtheorem{pr}[tm]{Proposition}

\makeatletter
\newcommand{\subscripts}[3]{%
  \@mathmeasure\z@\displaystyle{#2}%
  \global\setbox\@ne\vbox to\ht\z@{}\dp\@ne\dp\z@
  \setbox\tw@\box\@ne
  \@mathmeasure4\displaystyle{\copy\tw@_{#1}}%
  \@mathmeasure6\displaystyle{{#2}_{#3}}%
  \dimen@-\wd6 \advance\dimen@\wd4 \advance\dimen@\wd\z@
  \hbox to\dimen@{}\mathop{\kern-\dimen@\box4\box6}%
}
\makeatother

\newcommand{\T}{\mathcal{T}}
\newcommand{\g}{\frak{g}}
\newcommand{\nn}{\nonumber}

\newcommand{\h}{\mathrm{H}}

\newcommand{\ve}{\varepsilon}
\newcommand{\dis}{\displaystyle}
\newcommand{\Ker}{\mathrm{Ker}\,}

\newcommand{\del}{\partial}
\newcommand{\ol}{\overline}
\newcommand{\la}{\langle}
\newcommand{\La}{\langle\!\langle}
\newcommand{\ra}{\rangle}
\newcommand{\K}{\mathcal{K}}
\newcommand{\Ra}{\rangle\!\rangle}
\newcommand{\LA}{\longrightarrow}
\newcommand{\Span}{\mathrm{span}_{\mathbb{R}}}
\newcommand{\G}{G_{\Gamma}}

\newcommand{\Hom}{\mathrm{Hom}}

\newcommand{\Log}{\text{\rm{\,log\,}}}
\newcommand{\Exp}{\text{\rm{exp\,}}}

\newcommand{\AC}{\mathrm{AC}_0([0, 1]; \g^{(1)})}
\newcommand{\Conti}{C_0([0, 1]; \g^{(1)})}
\newcommand{\asy}{\rho_{\mathbb{R}}(\gamma_p)}
\newcommand{\Map}{\mathrm{Map}_0([0, 1]; \g^{(1)})}
\newcommand{\El}{L^\infty([0, 1]; \g^{(1)})}
\newcommand{\ContiG}{C_{\bm{1}_G}([0, 1]; G)}

\allowdisplaybreaks[3]

\makeatletter
 
 \@addtoreset{equation}{section}
\makeatother

\begin{document}
\title[Laws of the iterated logarithm on covering graphs]{
Laws of the iterated logarithm on covering graphs \\
with groups of polynomial volume growth
}
\author[Ryuya Namba]{Ryuya Namba}
\date{\today}
\address{Department of Mathematical Sciences,
College of Science and Engineering,
Ritsumeikan University, 1-1-1, Noji-higashi, Kusatsu
525-8577, Japan}
\email{{\tt{rnamba@fc.ritsumei.ac.jp}}}

\maketitle 
%
%
\begin{abstract}
Moderate deviation principles (MDPs)
for random walks on covering graphs
with groups of polynomial volume growth are discussed
in a geometric point of view. They deal with
any intermediate spatial scalings between those of laws of large numbers 
and those of central limit theorems. 
The corresponding rate functions are given by 
quadratic forms determined by   
the Albanese metric associated with the given random walks. 
We apply MDPs to establish
laws of the iterated logarithm on the covering graphs by
characterizing the set of all limit points of the normalized random walks. 
\end{abstract}


\section{Introduction and main results}

Let $K$ be a positive integer and $N=(N, \cdot)$ a finitely generated group with a generating set 
$\mathcal{S}=\{\gamma_1^{\pm1}, \gamma_2^{\pm1}, \dots, \gamma_K^{\pm 1}\}$.
We assume that the unit element 
$\bm{1}_N$ is not included in $\mathcal{S}$ for simplicity.  
A finitely generated group $N$ is said to be of {\it polynomial volume growth}
if there exist constants 
$C>0$ and $A \in \mathbb{N}$ such that 
$$
\#\{\gamma_{i_1}^{\ve_1} \cdot \gamma_{i_2}^{\ve_2} \cdot 
\cdots \cdot \gamma_{i_n}^{\ve_n} \, | \, 
\ve_k=\pm 1, \, i_k=1, 2, \dots, K, \, k=1, 2, \dots, n\} \leq Cn^{A}, \qquad n \in \mathbb{N}.
$$
The minimum of such integers $A$ is called the {\it growth rate} of $N$. 

Let $X=(V, E)$ be a covering graph 
of a finite graph  with a finitely generated 
covering transformation group of polynomial volume growth $N$. 
Here $V$ is the set of all vertices and $E$ is the set of all oriented edges. 
Gromov's celebrated theorem (see \cite{Gromov}) asserts that 
such $N$ is a finite extension of a torsion free nilpotent subgroup. 
Let $\Gamma$ be a torsion free nilpotent subgroup of $N$. 
Then $X$ is a covering graph of the finite graph $X_0=\Gamma \backslash X$ whose 
covering transformation group is $\Gamma$. 
Therefore, in the present paper, 
$X$ may be regarded as a {\it nilpotent covering graph} of a finite graph $X_0$
with a covering transformation group $\Gamma$
without loss of generality. A typical example of such groups is 
the {\it 3-dimensional Heisenberg group} defined by 
$$
\mathbb{H}^3(\mathbb{Z}):=\Bigg\{ \begin{pmatrix}
1 & x & z \\
0 & 1 & y \\
0 & 0 & 1 \end{pmatrix}
\Bigg| \, x, y, z \in \mathbb{Z}\Bigg\}.
$$
It is in some sense the simplest nilpotent group and 
its growth rate is known to be $A=4$. 
If $\Gamma$ is in particular an abelian group, then
$X$ is called a {\it crystal lattice}, which includes
classical periodic graphs such as square lattices, 
triangular lattices, hexagonal lattices and so on.  
We emphasize that crystal lattices have periodic global structures, 
while they have inhomogeneous local structures. 

On the other hand, it is common to consider {\it Cayley graphs}
of finitely generated groups in the study of random walks on groups. 
Generally speaking, as a finitely generated group is apart from abelian, 
the analysis on the corresponding Cayley graph becomes more complex.
However, we can say that nilpotent groups are ``almost abelian'' so that 
random walks on nilpotent Cayley graphs have been studied extensively 
(see e.g., \cite{A2} for the precise analysis of their long time asymptotics). 
We here note that a nilpotent Cayley graph is regarded as a nilpotent covering graph 
of a finite bouquet graph with a nilpotent covering transformation group. 
Therefore, it is natural to deal with nilpotent covering graphs as hybrid models of both
crystal lattices and nilpotent Cayley graphs. 

We are interested in investigating long time asymptotics of random walks 
on a nilpotent covering graph $X$.
However, it is difficult to discuss such problems on $X$ directly since
$X$ is a discrete model on which any spatial scalings cannot defined. 
To overcome the difficulties, 
we employ Malc\'ev's theorem (cf.~\cite{Malcev}), which
asserts that there exists a connected and simply connected 
nilpotent Lie group $(G, \cdot )=\G$ 
such that $\Gamma$ is isomorphic to a cocompact lattice in $G$. 
We then come to discuss limit theorems on $X$
through a $\Gamma$-equivariant map $\Phi : X \to G$, that is, 
$\Phi(\gamma x)=\gamma \cdot \Phi(x)$ for $\gamma \in \Gamma$
and $x \in X$. Such a map $\Phi$ is called a {\it periodic realization} of $X$. 
For example, there are several studies of
{\it central limit theorems} (CLT) for random walks on $X$ 
in a point of view of {\it discrete geometric analysis} initiated by Sunada 
(see e.g.,~\cite{S} and \cite{KS06} for overviews of the theory). 
Roughly speaking, it is the theory of discrete harmonic analysis on periodic graphs. 
As for the case where $X$ is a crystal lattice, the authors in \cite{IKK} showed
{\it functional} CLTs (i.e.,\,Donsker-type invariance principles)
for non-symmetric random walks on $X$ by using 
the notion of {\it modified harmonic realizations} 
$\Phi_0 : X \to \Gamma \otimes \mathbb{R} \cong \mathbb{R}^d$
introduced in \cite{KS06}. 
See also \cite{Kotani} for CLTs for magnetic transition operators on $X$. 
As a generalization to nilpotent cases, 
Ishiwata, Kawabi and Namba proved functional CLTs for non-symmetric random walks 
on a nilpotent covering graph $X$
in \cite{IKN-1, IKN-2}. 
They extended the modified harmonicity to nilpotent settings, which
also plays crucial roles in the proof of the CLTs. 

On the other hand, 
Kotani and Sunada showed in \cite{KS06} that
random walks on crystal lattices satisfy
{\it large deviation principles} (LDPs) 
by employing discrete geometric analysis developed by themselves.
We also refer to \cite{Kotani contemp} for related results on LDPs on crystal lattices. 
In nilpotent cases, Baldi and Caramellino \cite{BC} showed
LDPs for random walks given by the product of i.i.d.\,random variables  
on a nilpotent Lie group $N$. 
Note that the random walks on $N$ is regarded as the ones
on nilpotent covering graphs of finite bouquet graphs
whose covering transformation group is a proper nilpotent lattice in $N$. 
Afterwards, LDPs on nilpotent covering graphs were discussed in \cite{Tanaka}
as generalizations of both \cite{KS06} and \cite{BC}, 
which basically motivates a part of our study.  
Let us consider an random walk $\{w_n\}_{n=0}^\infty$ with values in 
a nilpotent covering graph $X$, 
which is not always symmetric. Write $\mathbb{P}_x$ for its probability law 
when it starts at $x \in X$. 
A periodic realization $\Phi : X \to G$ then induces 
a $G$-valued random walk $\{\xi_n=\Phi(w_n)\}_{n=0}^\infty$.
Note that the realization $\Phi$ is not necessarily harmonic, 
though the harmonicity plays a key role in
proving CLTs as we mentioned above.

Since $G$ is nilpotent, the corresponding Lie algebra $\g$ may have the decomposition
$\g=\g^{(1)} \oplus \g^{(2)} \oplus \cdots \oplus \g^{(r)}$
for some $r \in \mathbb{N}$, satisfying that  $[\g^{(i)}, \g^{(j)}] \subset \g^{(i+j)}$, 
$i+j=2, 3, \dots, r$, and $\g^{(1)}$ generates $\g$ (see Section 2.1). 
We occasionally call $\g^{(1)}$ the generating part of the Lie algebra $\g$. 
We note that $\g^{(1)} \sim G/[G, G]$, which means that $\g^{(1)}$ is regarded as 
the abelianization of $G$.

Let us go back to the random walk $\{\xi_n=\Phi(w_n)\}_{n=0}^\infty$ on $G$.
Since $\log(\xi_n)$ takes its values in $\g$, it can be uniquely decomposed as
$$
\log(\xi_n)=\log(\xi_n)|_{\g^{(1)}}+\log(\xi_n)|_{\g^{(2)}}+\cdots+\log(\xi_n)|_{\g^{(r)}},
$$ 
where $\log(\xi_n)|_{\g^{(i)}} \in \g^{(i)}$ for $i=1, 2, \dots, r$. 
Here, $\log : G \to \g$ means
the inverse map of the usual exponential map $\exp : \g \to G$. 
We set
$\Xi_n:=\log(\xi_n)\big|_{\g^{(1)}}$ for $n=0, 1, 2, \dots$,
which gives a random walk on $\g^{(1)}$. 

We should note that a {\it law of large numbers} (LLN) holds for 
$\{\Xi_n\}_{n=0}^\infty$ (cf. \cite[Section 3.4]{IKN-1}): we obtain
\begin{equation}\label{LLN}
\lim_{n \to \infty}\frac{1}{n}\Xi_n=\rho_{\mathbb{R}}(\gamma_p) \in \g^{(1)} \quad \mathbb{P}_x
\text{-a.s.},
\end{equation}
where the quantity $\asy$ is called the {\it asymptotic direction} of the given random walk.
The quantity $\asy$ is regarded as the ``mean'' of the random walk $\{\Xi_n\}_{n=0}^\infty$. 
Indeed, if the random walk $\{w_n\}_{n=0}^\infty$ is symmetric, then it holds that 
$\asy=\bm{0}$. 
See Section 2.2 for more details on the quantities $\gamma_p$ and $\rho_{\mathbb{R}}(\gamma_p)$.  
Let $G_\infty=(G, *)$ be the {\it limit group} of $G$  and 
$\varphi : G \to G_\infty$ a canonical diffeomorphism (see Section 2.1). 
We denote by $(\tau_\ve)_{\ve \geq 0}$ the one-parameter group of 
{\it dilations} on $G_\infty$, which behaves like scalar multiplications on $G_\infty$. 
Then the LDP established in \cite{Tanaka} is stated as follows:

\begin{pr}[{\bf cf.~\cite[Theorem 1.1]{Tanaka}}]\label{LDP on G}
An LDP holds for the sequence of $G_\infty$-valued random variables 
$\{\tau_{1/n}\varphi(\xi_n)\}_{n=0}^\infty$ with the speed rate $n$ 
and with some good rate function $I_\infty : G_\infty \to [0, \infty]$. 
Namely, it holds that
$$
\begin{aligned}
-\inf_{y \in A^\circ}I_\infty(y) &\leq \liminf_{n \to \infty}\frac{1}{n}
\log\mathbb{P}_x\Big(\tau_{1/n}\big(\varphi(\xi_n)\big) \in A\Big)\\
&\leq  \limsup_{n \to \infty}\frac{1}{n}\log\mathbb{P}_x\Big(\tau_{1/n}\big(\varphi(\xi_n)\big) \in A\Big)
\leq -\inf_{y \in \ol{A}}I_\infty(y)
\end{aligned}
$$
for a Borel set $A \subset G_\infty$, where 
$A^\circ$ and $\ol{A}$ are the interior and the closure of $A$, respectively.  
\end{pr}
\noindent
By taking (\ref{LLN}) into account, we can say that Proposition \ref{LDP on G} deals with
exactly how fast
the probability
$$
\mathbb{P}_x\Big(d_{\mathrm{Fin}}\Big( \tau_{1/n}\big(\varphi(\xi_n)\big), \exp\big(\asy) \Big)>\delta\Big),
\qquad \delta>0
$$
exponentially decays as $n \to \infty$, 
where $d_{\mathrm{Fin}}$ is a canonical left invariant metric on $G_\infty$ called the {\it Finsler metric}, 
whose definition will be stated in Section 3.

The first aim of the present paper is to establish 
{\it moderate deviation principles} (MDPs) for a random walk
$(\Omega_x(X), \mathbb{P}_x, \{w_n\}_{n=0}^\infty)$
on a $\Gamma$-nilpotent covering graph $X$, where $\Omega_x(X)$ is the set 
of all paths in $X$ starting at $x \in V$. 
Formally there seems to be no difference between MDPs and LDPs.  
MDPs deal with what happens at any intermediate speed rates
between $n$ of LLN-type and $\sqrt{n}$ of CLT-type, while
LDPs usually deal with the speed rate $n$ of LLN-type. 
It is known that MDPs cannot be obtained neither LDPs  
nor CLTs. Baldi and Caramellino \cite{BC} also showed
MDPs for the product of i.i.d.\,and {\it centered} random variables  
on a nilpotent Lie group $N$. However, MDPs on covering graphs 
have not been obtained because random walks 
on covering graphs are not always i.i.d. 
Our main concern is to know how fast the probability
$$
\mathbb{P}_x\Big(d_{\mathrm{Fin}}\Big( \varphi(\xi_n), \exp\big(n\asy) \Big)>\delta a_n\Big),
\qquad \delta>0
$$
decays exponentially, where $\{a_n\}_{n=1}^\infty$ is an arbitrary monotone increasing 
sequence of positive real numbers satisfying 
\begin{equation}\label{scale}
\lim_{n \to \infty}\frac{a_n}{\sqrt{n}}=+\infty \quad \text{and}
\quad \lim_{n \to \infty}\frac{a_n}{n}=0. 
\end{equation}
 
We introduce two sequences of $\g^{(1)}$-valued random variables 
$\{W_n\}_{n=1}^\infty$ and $\{\ol{W}_n\}_{n=1}^\infty$ by setting
$$
W_n(c):=\Xi_{n}(c) - \Xi_{n-1}(c), \qquad \ol{W}_n(c):=W_n(c) - \asy, 
\qquad n \in \mathbb{N}, \, c \in \Omega_x(X),
$$
respectively.
Let $\AC$ be the set of all absolutely continuous paths $h : [0, 1] \to \g^{(1)}$ with $h(0)=0$.
We then define an $\AC$-valued random variable $Z^{(n)} : \Omega_x(X) \to \AC$ by
$$
Z_t^{(n)}(c):=\frac{1}{a_n}\ol{W}_1(c)+\frac{1}{a_n}\ol{W}_2(c)+\cdots + \frac{1}{a_n}\ol{W}_{[nt]}(c)
+\frac{1}{a_n}(nt-[nt])\ol{W}_{[nt]+1}(c)
$$
for $n \in \mathbb{N}$, $0 \leq t \leq 1$ and $c \in \Omega_x(X)$. 
Let $g_0$ be the {\it Albanese metric} on $\g^{(1)}$ associated with the given random walk
(see Section 2.2 for details). We put $\Sigma^{-1}=(\la X_i^{(1)}, X_j^{(1)}\ra_{g_0})_{i, j=1}^{d_1}$,
where $\{X_1^{(1)}, X_2^{(1)}, \dots, X_{d_1}^{(1)}\}$ is a fixed basis of $\g^{(1)}$. 
In order to establish our first main result,
We need to show an MDP for the sequence $\{Z^{(n)}\}_{n=1}^\infty$. 

\begin{pr}[{\bf MDP on a path space}] \label{MDP on a path space}
The sequence of $\AC$-valued random variables $\{Z^{(n)}\}_{n=1}^\infty$
satisfies an MDP with the speed rate $a_n^2/n$. 
Namely,  it holds that  
$$
\begin{aligned}
-\inf_{h \in A^\circ}I'(h) &\leq \liminf_{n \to \infty}\frac{n}{a_n^2}\log\mathbb{P}_x(Z^{(n)} \in A)\\
&\leq  \limsup_{n \to \infty}\frac{n}{a_n^2}\log\mathbb{P}_x(Z^{(n)} \in A)
\leq -\inf_{h \in \ol{A}}I'(h)
\end{aligned}
$$
for any Borel set $A \subset \AC$, where $I' : \AC \to [0, \infty]$ is a good rate function defined by
\begin{equation}\label{rate function path space}
I'(h):=\begin{cases}
\dis \int_0^1 \alpha^*\big(\dot{h}(t)\big) \, dt & \text{if }h \in \AC, \\
+\infty & \text{otherwise}.
\end{cases}
\end{equation}
Here, $\alpha^* : \g^{(1)} \to \mathbb{R}$ is the quadratic form given by $\alpha^*(\chi):=\frac{1}{2}\la \Sigma^{-1}\chi, \chi\ra$ for $\chi \in \g^{(1)}$. 
\end{pr}
\noindent
We should mention that MDPs on path spaces have been investigated in various settings. 
Mogulskii first showed that trajectories of the sums of $\mathbb{R}^d$-valued i.i.d.\,random variables
satisfy MDPs in \cite{Mog}. 
Afterwards, the authors generalized MDPs above to the case of
trajectories of the sums of i.i.d.\,random vectors with values in Hausdorff topological vector spaces
in \cite{BM}. 
We also refer to 
\cite{HL} for extensions of such MDPs to infinite-dimensional cases.

We now put 
$\ol{\xi}_n:=\xi_n \cdot \exp(-n \asy)$ for $n \in \mathbb{N}$. 
We again note that our underlying random walk $\{w_n\}_{n=0}^\infty $ 
is not always symmetric. 
If $\{w_n\}_{n=0}^\infty $ is symmetric, then we have $\ol{\xi}_n:=\xi_n $ for $n \in \mathbb{N}$.
By combining Proposition \ref{MDP on a path space} with several well-known facts 
in the theory of large deviations, 
we obtain the following. 

\begin{tm}[{\bf MDP on $X$}] \label{MDP on G}
The sequence of $G_\infty$-valued random variables 
$\{\tau_{1/a_n}\varphi(\ol{\xi}_n)\}_{n=0}^\infty$ 
satisfies an MDP
with the speed rate $a_n^2/n$ and a good rate function $I_\infty : G_\infty \to [0, \infty]$.
Namely, it holds that
$$
\begin{aligned}
-\inf_{g \in A^\circ}I_\infty(g) &\leq \liminf_{n \to \infty}\frac{n}{a_n^2}
\log\mathbb{P}_x\Big(\tau_{1/a_n}\big(\varphi(\ol{\xi}_n)\big) \in A\Big)\\
&\leq  \limsup_{n \to \infty}\frac{n}{a_n^2}\log\mathbb{P}_x
\Big(\tau_{1/a_n}\big(\varphi(\ol{\xi}_n)\big) \in A\Big)
\leq -\inf_{g \in \ol{A}}I_\infty(g)
\end{aligned}
$$
for any Borel set $A \subset G_\infty$.
\end{tm}
\noindent
We emphasize that the rate function $I_\infty$ can be 
written in terms of $I'$ defined by (\ref{rate function path space})
and it does not depend on 
the choice of the sequence $\{a_n\}_{n=1}^\infty \subset \mathbb{R}$ satisfying (\ref{scale}).
The precise definition of $I_\infty$ will be clearly stated in Proposition \ref{MDP pre}.

We here give a few remarks on Theorem \ref{MDP on G}. 
We first note that Theorem $\ref{MDP on G}$ deals with the sequence 
$\{\ol{\xi}_n\}_{n=0}^\infty$ instead of $\{\xi_n\}_{n=0}^\infty$, 
while Proposition \ref{LDP on G}
directly deals with the sequence $\{\xi_n\}_{n=0}^\infty$. 
In fact, it is common to impose the ``centered'' 
condition to random variables in MDP frameworks.
That plays a crucial role to show that the first order terms of the logarithmic 
moment generating function of $Z^{(n)}$ converges to zero as $n \to \infty$ 
under suitable spatial scalings in Lemma \ref{MDP I_j} below.
The second remark comes from geometric perspectives. 
It turns out in \cite{Tanaka} that Proposition \ref{LDP on G} is closely related to
the pointed Gromov--Hausdorff limit of nilpotent covering graphs 
with suitably normalized graph distance (see also \cite{Pansu}). 
Furthermore, a characterization of the effective domain 
of the rate function in Proposition \ref{LDP on G}
is given in terms of the Carnot--Carath\'eodory metric 
$d_{\mathrm{CC}}$ (cf.~\cite[Theorem 1.2]{Tanaka}). 
Therefore, one may wonder if there are or are not relations 
between Theorem \ref{MDP on G} 
and the geometry of underlying graphs. 
We have obtained an explicit representation of 
the rate function $I_\infty$ in Theorem \ref{MDP on G}
in terms of the Albanese metric $g_0$ on $\g^{(1)}$. 
We provide a first result of MDPs with an explicit rate function
in a geometric sense.  
As a further problem, to investigate relations between 
Theorem \ref{MDP on G} and some convergences of nilpotent 
covering graphs might be interesting.

From now on, we consider a certain sequence $b_n=\sqrt{n \log \log n}$ for $n \in \mathbb{N}$. 
It is needless to say that $\{b_n\}_{n=1}^\infty$ is a typical example which satisfies (\ref{scale}), 
so that Theorem \ref{MDP on G} implies that an MDP holds for 
$\{\tau_{1/b_n} (\varphi(\ol{\xi}_n))\}_{n=1}^\infty$ with the speed rate $b_n^2/n=\log\log n$.
On the other hand, it is known that the sequence $\{b_n\}_{n=1}^\infty$
has an explicit interaction with {\it laws of the iterated logarithm} (LILs, in short). 
Let us review a few known results on LILs in nilpotent settings. 
As an early work, there is a paper \cite{CR} in which 
an LIL for random walks on the 3-dimensional Heisenberg group $N=\mathbb{H}^3(\mathbb{R})$ 
was established under some moment conditions. 
We also refer to \cite{Neu} for related limit theorems on $\mathbb{H}^3(\mathbb{R})$, as well as LILs.
Baldi and Caramellino mentioned in \cite{BC} that LILs for random walks 
given by the product of i.i.d.\,random variables on nilpotent Lie groups
might be established by applying MDPs obtained by themselves. 
Nevertheless, the complete 
proof was not given in the paper. 
Afterwards, the authors gave a complete proof of LILs on nilpotent Lie groups in \cite{CV},
by characterizing the set of all limit points of random walks scaled by $\{b_n\}_{n=1}^\infty$.
Such kinds of LILs are known to {\it LILs of functional type}. 

In spite of these developments, LILs for random walks on 
nilpotent covering graphs have not been studied, even in cases of crystal lattices. 
In fact, it is difficult to establish them since
random walks on covering graphs 
are not always i.i.d., while to consider i.i.d.\,cases plays 
essential roles to prove LILs in \cite{CR, BC, CV}.
Under these circumstances,
the second purpose of the present paper is to prove LILs of functional type
for random walks on a
nilpotent covering graph $X$, by identifying the set of all limit points of 
the sequence $\{\tau_{1/b_n}(\varphi(\ol{\xi}_n))\}_{n=0}^\infty$. 
The precise statement is as follows:

\begin{tm}[{{\bf LIL on $X$}}]\label{LIL}
Let $b_n=\sqrt{n \log \log n}$ for $n \in \mathbb{N}$. 
We denote by $\widehat{\mathcal{K}}$ the set of all $\mathbb{P}_x$-a.s.\,limit points of
the sequence $\{\tau_{1/b_n}(\varphi(\ol{\xi}_n))\}_{n=0}^\infty$ of $G_\infty$-valued random variables. 
Then we obtain
$$
\widehat{\mathcal{K}}=\{g \in G_\infty \, | \, I_\infty(g) \leq 1\},
$$
where $I_\infty : G_\infty \to [0, \infty]$ is the rate function given in Theorem \ref{MDP on G}. 
\end{tm}
Towards the proof, we first need to show that 
the sequence of $\AC$-valued random variables $\{Z^{(n)}\}_{n=1}^\infty$
is relatively compact in $\Conti$, the set of all continuous paths 
taking values in $\g^{(1)}$ and zero at the origin (Lemma \ref{LIL pre}). 
Let $\mathcal{K}$ be the all $\mathbb{P}_x$-a.s.\,limit points of the sequence $\{Z^{(n)}\}_{n=1}^\infty$.
We next describe a explicit representation of the set $\mathcal{K}$ in terms of the rate function 
$I' : \AC \to [0, \infty]$ (Lemma \ref{LIL-pre2}) by applying Proposition \ref{MDP on a path space}. 
Finally, a proof of Theorem \ref{LIL} is given 
by combining these two lemmas with Theorem \ref{MDP on G}. 

\begin{re}
There are several contexts in which a sequence 
$\widetilde{b}_n:=\sqrt{2n \log \log n}$ for $n \in \mathbb{N}$ is used
in considering LILs (see e.g., \cite{BC, HL}). 
If we choose $\{\widetilde{b}_n\}_{n=1}^\infty$ instead of $\{b_n\}_{n=1}^\infty$, then
we can show that 
$$
\widehat{\mathcal{K}}'=\{g \in G_\infty \, | \, 2I_\infty(g) \leq 1\}
$$
appears as the limit set of random walks normalized by 
$\{\widetilde{b}_n\}_{n=0}^\infty$ insetad of $\widehat{\mathcal{K}}$. 
\end{re}


\section{Notations}


\subsection{Nilpotent Lie groups and their limit groups}
Let $G$ be a connected and simply connected nilpotent Lie group 
and $\g$ the corresponding Lie algebra. 
We now set
$\frak{n}_1:=\frak{g}$ and $\frak{n}_{k+1}:=[\frak{g}, \frak{n}_k]$ for $k \in \mathbb{N}$.
Since $\g$ is nilpotent, 
there exists an integer $r \in \mathbb{N}$ such that
$\frak{g}=\frak{n}_1 \supset \dots 
\supset \frak{n}_r \supsetneq \frak{n}_{r+1}=\{\bm{0}_{\frak{g}}\}$,
which is called the {\it step number} of $\g$ or $G$. 
We define the subspace $\frak{g}^{(k)}$ of $\frak{g}$ by
$\frak{n}_k=\frak{g}^{(k)} \oplus \frak{n}_{k+1}$ for $k=1, 2, \dots, r$.
Then the Lie algebra $\frak{g}$ is decomposed as 
$\frak{g}=\frak{g}^{(1)} \oplus \g^{(2)} \oplus \dots \oplus \frak{g}^{(r)}$ 
and each  $Z \in \frak{g}$ is uniquely written as $Z=Z^{(1)} + Z^{(2)}+\dots + Z^{(r)}$, 
where 
$Z^{(k)} \in \frak{g}^{(k)}$ for $k=1, 2, \dots, r$.
We define a map $T_\ve : \g \LA \g$ by
$$
T_\ve(Z):=\ve Z^{(1)} + \ve^2 Z^{(2)} + \dots + \ve^r Z^{(r)}, \qquad \ve \geq 0, \, Z \in \frak{g}
$$
and also define a Lie bracket product $[\![ \cdot, \cdot ]\!]$ on $\frak{g}$ by
$$
[\![ Z_1, Z_2 ] \!]:=\lim_{\ve \searrow 0} T_\ve
 \big[T_{1/\ve}(Z_1), T_{1/\ve}(Z_2)\big],
\qquad Z_1, Z_2 \in \frak{g}.
$$ 
We then introduce a map $\tau_\ve : G \LA G$, called the {\it dilation operator} on $G$, defined by
$$
\tau_\ve (g):=\Exp \big(T_\ve \big( \Log(g)\big)\big), \qquad \ve \geq 0, \, g \in G,
$$
which gives scalar multiplications on $G$. 
We note that $\tau_\ve$ may not be a group homomorphism, though 
it is a diffeomorphism on $G$. 
By the dilation map $\tau_\ve$, a Lie group product $*$ on $G$ 
is defined by
$$
g * h :=\lim_{\ve \searrow 0} 
\tau_{\ve}\big( \tau_{1/\ve}(g) \cdot \tau_{1/\ve}(h)\big), \qquad g, h \in G.
$$
The Lie group $G_\infty=(G, *)$ is called the {\it limit group} of $(G, \cdot)$. 
The corresponding Lie algebra $\g_\infty=(\g, [\![\cdot, \cdot]\!])$
is {\it stratified} in the sense that it is decomposed as
$\g=\bigoplus_{k=1}^r \frak{g}^{(k)}$ satisfying that
$[\![\frak{g}^{(k)}, \frak{g}^{(\ell)}]\!] \subset \g^{(i+j)}$ unless 
$i+j>r$
and the subspace $\frak{g}^{(1)}$ generates $\frak{g}$. 
It is known that the exponential map $\exp_\infty : \g_\infty \to G_\infty$
coincides with $\exp : \g \to G$ as a map. 
Moreover, we see that $G_\infty$ is diffeomorphic to $G$
in the following manner. 
The construction of $\g_\infty$ immediately induces a canonical linear map
$\iota : \g \to \g_\infty$. We then define a map $\varphi : G \to G_\infty$ by
$\varphi:=\exp_\infty \circ \iota \circ \exp^{-1}$. 
This map $\varphi$ does gives a canonical diffeomorphism between $G$ and $G_\infty$. 
Let $\{X_1^{(1)}, X_2^{(1)}, \dots, X_{d_1}^{(1)}\}$ be a fixed basis of the generating part $\g^{(1)}$. 
We define a norm $\|\cdot\|_{\g^{(1)}}$ on $\g^{(1)}$ by 
$$
\|Z\|_{\g^{(1)}}:=\Big(\sum_{i=1}^{d_1} c_i^2\Big)^{1/2} \quad \text{when } Z=\sum_{i=1}^{d_1}c_iX_i^{(1)} \in \g^{(1)}. 
$$
In what follows, a basis $\{X_1^{(1)}, X_2^{(1)}, \dots, X_{d_1}^{(1)}\}$ of $\g^{(1)}$
is supposed to be fixed. By introducing an arbitrary norm 
$\|\cdot\|$ on $\g^{(2)} \oplus \g^{(3)} \oplus \cdots \oplus \g^{(r)}$, 
we define a norm $\|\cdot\|_{\g}$ on $\g$ by $\|Z\|_{\g}:=\|Z_1\|_{\g^{(1)}}+\|Z_2\|$, 
where $Z=Z_1+Z_2$ with $Z_1 \in \g^{(1)}$ and 
$Z_2 \in \g^{(2)} \oplus \g^{(3)} \oplus \cdots \oplus \g^{(r)}$. 

For more details on nilpotent Lie groups, see e.g., \cite{Goodman, VSC}. 
The two different ways to construct $G_\infty$ are discussed in \cite{IKN-1, Tanaka}. 
However, we emphasize that they are alternative ones so that 
there is no difference whichever we choose.

\subsection{Random walks and discrete geometric analysis}

Let $X=(V, E)$ be a $\Gamma$-nilpotent covering graph of a finite graph $X_0=(V_0, E_0)$.  
For an edge $e \in E$, 
we denote by $o(e)$, $t(e)$ and $\ol{e}$ 
the origin, the terminus and the inverse edge of $e$, respectively. 
Set $E_x=\{e \in E \, | \, o(e)=x\}$ for $x \in V.$
A {\it path} $c$ in $X$ of length $n$ is a sequence $c=(e_1, e_2, \dots, e_n)$ of $n$ edges 
$e_1, e_2, \dots, e_n \in E$ with $o(e_{i+1})=t(e_i)$ for $i=1, 2, \dots, n-1$. 
Let $\Omega_{x, n}(X)$ be
the set of all paths in $X$ of length $n \in \mathbb{N} \cup \{\infty\}$ starting from $x \in V$. 
Put $\Omega_x(X)=\Omega_{x, \infty}(X)$ for simplicity. 

Let $p : E_0 \to (0, 1]$ be a transition probability. Namely, it satisfies 
$\sum_{e \in (E_0)_x}p(e)=1$ for $x \in V_0$  and $p(e)>0$ for $e \in E_0$.
Then the transition probability yields 
a $X_0$-valued time-homogeneous Markov chain 
$(\Omega_x(X_0), \mathbb{P}_x, \{w_n\}_{n=0}^\infty)$, where $\mathbb{P}_x$ is 
the probability measure on $\Omega_x(X_0)$ induced from  $p$
and $w_n(c):=o( e_{n+1})$ for $n \in \mathbb{N} \cup\{0\}$
and $c=(e_1, e_2, \dots, e_n, \dots) \in \Omega_x(X_0)$. 
We find a {\it normalized invariant measure} $m : V_0 \to (0, 1]$
by applying the Perron--Frobenius theorem. 
We put $\widetilde{m}(e):=p(e)m\big(o(e)\big)$ for $e \in E_0$.
A random walk satisfying $\widetilde{m}(e)=\widetilde{m}(\ol{e})$ for $e \in E_0$
is said to be ($m$-){\it symmetric}. 
A random walk on $X$ is defined by a $\Gamma$-invariant lift of the 
random walk on $X_0$. Namely, the transition probability, say also $p : E \to (0, 1]$,
satisfies $p(\gamma e)=p(e)$ for $\gamma \in \Gamma$ and $e \in E$. 
We also write $\mathbb{P}_x$ the probability measure on $\Omega_x(X)$ induced by $p : E \to (0, 1]$.
If $c=(e_1, e_2, \dots, e_n) \in \Omega_{x, n}(X)$, then we put $p(c)=p(e_1)p(e_2) \cdots p(e_n)$. 
For more basics on random walks on graphs or groups, we refer to \cite{Woess}. 

We define the 0-chain group and 1-chain group by
$$
\begin{aligned}
C_0(X_0, \mathbb{R})&:=\Big\{ \sum_{x \in V_0}a_x x \, \Big| \, a_x \in \mathbb{R}\Big\}, &
C_1(X_0, \mathbb{R})&:=\Big\{ \sum_{e \in E_0}a_e e \, \Big| \, a_e \in \mathbb{R}, \, \ol{e}=-e\Big\}, 
\end{aligned}
$$
respectively. 
The boundary operator $\del : C_1(X_0, \mathbb{R}) \LA C_0(X_0, \mathbb{R})$ 
is defined by $\del(e)=t(e)-o(e)$ for $e \in E_0$. 
Then, the first homology group $\h_1(X_0, \mathbb{R})$ is defined by 
$\Ker(\del) \subset C_1(X_0, \mathbb{R})$. 
We define the {\it homological direction} of $X_0$  by
$$
\gamma_p:=\sum_{e \in E_0}\widetilde{m}(e)e \in \h_1(X_0, \mathbb{R}),
$$
which gives a homological drift of the given random walk. 
It is easily verified that a random walk on $X_0$ is ($m$-)symmetric if and only if $\gamma_p=0$. 
By employing the discrete analogue of Hodge--Kodaira theorem
(cf.~\cite[Lemma 5.2]{KS06}), we equip
the first cohomology group $\h^1(X_0, \mathbb{R}):=\big(\h_1(X_0, \mathbb{R})\big)^*$ 
with the inner product
$$
\La \omega_1, \omega_2 \Ra_p
:=\sum_{e \in E_0}\widetilde{m}(e)\omega_1(e)\omega_2(e)
 - \la \gamma_p, \omega_1 \ra \la \gamma_p, \omega_2 \ra,
 \qquad \omega_1, \omega_2 \in \h^1(X_0, \mathbb{R})
$$
associated with the transition probability $p$. 
Let $\rho_{\mathbb{R}} : \h_1(X_0, \mathbb{R}) \to \g^{(1)}$ be
the canonical surjective linear map
induced by the canonical surjective homomorphism $\rho : \pi_1(X_0) \to \Gamma$, 
where $\pi_1(X_0)$ is the fundamental group of $X_0$. 
We call the quantity $\rho_{\mathbb{R}}(\gamma_p) \in \g^{(1)}$
the {\it asymptotic direction} of $X_0$. 
We emphasize that $\asy$ has already appeared in (\ref{LLN}).  
Then, through the transpose ${}^t\rho_{\mathbb{R}}$, a flat metric $g_0$ on $\g^{(1)}$ is induced 
from $\La \cdot, \cdot \Ra_p$ as in the diagram below. 
$$
\xymatrix{ 
(\frak{g}^{(1)}, g_0)  \ar @{<->}[d]^{\mathrm{dual}}
 \ar @{<<-}[r]^{\rho_{\mathbb{R}}} \ar @{<->}[d]^{\mathrm{dual}} 
 & \h_1(X_0, \mathbb{R}) \ar @{<->}[d]^{\mathrm{dual}} &\\
\Hom(\frak{g}^{(1)}, \mathbb{R}) 
\ar @{^{(}->}[r]_{{}^t \rho_{\mathbb{R}}\qquad} & (\h^1(X_0, \mathbb{R}), \La \cdot, \cdot \Ra_p).
}
$$
We call the metric $g_0$ the {\it Albanese metric}.

\begin{re}
In the present paper, there are two adjectives ``symmetric'' and ``centered'' for random walks. 
Since both are very alike, some readers might confuse their precise meanings. 
Therefore, it is worth rementioning a relation between these two notions in terms of the quantities 
$\gamma_p$ and $\rho_{\mathbb{R}}(\gamma_p)$. 
\begin{quote}

$\bullet$ $\{w_n\}_{n=0}^\infty$ is symmetric if $\gamma_p=0$.

\noindent
$\bullet$ $\{w_n\}_{n=0}^\infty$ is non-symmetric if $\gamma_p \neq 0$.

\noindent
$\bullet$ $\{w_n\}_{n=0}^\infty$ is centered if $\asy=\bm{0}$.

\noindent
$\bullet$ $\{w_n\}_{n=0}^\infty$ is non-centered if $\asy \neq \bm{0}$.

\end{quote}
It is easily seen that $\gamma_p=0$ implies $\asy=\bm{0}$. However, 
the converse does not hold in general. Namely, there is a case where
$\gamma_p \neq 0$ but $\asy =\bm{0}$ {\rm(} see e.g., \cite[Section 6]{IKN-1}{\rm)}. 
This is why we need to distinguish these two notions explicitly. 
We emphasize that the main results in the present paper 
hold even in the non-centered cases.  
\end{re}

Let $\{X_1^{(1)}, X_2^{(1)}, \dots, X_{d_1}^{(1)}\}$ be a fixed basis of $\g^{(1)}$
and $\{\omega_1, \omega_2, \dots, \omega_{d_1}\}$ the corresponding basis of
the dual space
$\Hom(\g^{(1)}, \mathbb{R}) \hookrightarrow (\h^1(X_0, \mathbb{R}), \La \cdot, \cdot \Ra_p)$.
By definition, we easily see that 
the ${d_1} \times d_1$-matrix $\Sigma:=\big( \La \omega_i, \omega_j \Ra_p\big)_{i, j=1}^{d_1}$
is invertible and its inverse matrix is given by 
$\Sigma^{-1}=\big( \la X_i^{(1)}, X_j^{(1)} \ra_{g_0} \big)_{i, j=1}^{d_1}$.
It is needless to say that $\Sigma$ is the identity (and so is  $\Sigma^{-1}$) if
we replace $\{X_1^{(1)}, X_2^{(1)}, \dots, X_{d_1}^{(1)}\}$ by an orthonormal basis
$\{V_1^{(1)}, V_2^{(1)}, \dots, V_{d_1}^{(1)}\}$ with respect to the Albanese metric $g_0$. 

We put
$\ol{\Xi}_n(c):=\Xi_n(c) - n\asy$ for $n=0, 1, 2, \dots$ and $c \in \Omega_x(X)$. 
We need to show the following lemma to know the long time behavior of 
the covariances of $\{\ol{\Xi}_n\}_{n=0}^\infty$, 
since $\{\ol{\Xi}_n\}_{n=0}^\infty$ is not i.i.d. in general (but is always independent by definition). 

\begin{lm}\label{convergence of cov}
For $1 \leq i, j \leq d_1$, we have
\begin{align}
&\lim_{N \to \infty}\frac{1}{N}
\mathbb{E}^{x}\Big[ 
 \ol{\Xi}_N\big|_{X_i^{(1)}} 
 \ol{\Xi}_N\big|_{X_j^{(1)}} \Big]
=\La \omega_i, \omega_j \Ra_p,
\end{align}
where $\mathbb{E}^x$ stands for the expectation with respect to $\mathbb{P}_x$. 
\end{lm}

\noindent
{\bf Proof.} Let $\Phi_0 : X \to G$ be the {\it modified harmonic realization}, that is, 
\begin{equation}\label{harmonicity}
\sum_{e \in E_x}p(e)\Big\{ \log \big(\Phi_0(t(e))\big)\big|_{\g^{(1)}}
- \log \big(\Phi_0(o(e))\big)\big|_{\g^{(1)}}\Big\}=\rho_{\mathbb{R}}(\gamma_p),
\qquad x \in V.
\end{equation}
The modified harmonic realization provides the most natural configuration of
vertices in the nilpotent Lie group $G$ in a geometric sense, which strongly relates 
to some probabilistic notions such as martingales. 
For more details of modified harmonic realizations, see \cite[Section 3.3]{IKN-1}. 
We then define a $\g^{(1)}$-valued random walk $\{\Xi^{(0)}_n\}_{n=0}^\infty$ by
$\Xi^{(0)}_n(c):=\log\big(\Phi_0(w_n(c)\big)\big|_{\g^{(1)}}$ 
and put 
$$
\ol{\Xi}^{(0)}_n(c):=\Xi^{(0)}_n(c) - n\asy
$$
for $n=0, 1, 2, \dots$ and $c \in \Omega_x(X)$. 
We also introduce a $\g^{(1)}$-{\it corrector} $\Psi : X \to \g^{(1)}$ of $\Phi$ by 
$$
\Psi(x):=\log \big(\Phi(x)\big)\big|_{\g^{(1)}} - \log \big(\Phi_0(x)\big)\big|_{\g^{(1)}},
\qquad x \in V,
$$
which plays a role in measuring the difference between $\Phi$ and $\Phi_0$. 
Note that the set $\{\Psi(x) \, | \, x \in V\}$ is finite due to 
$\Psi(\gamma x)=\Psi(x)$ for $\gamma \in \Gamma$ and $x \in V$. 
In particular, there exists a sufficiently large $M>0$ with
$\max_{x \in V} \|\Psi(x)\|_{\g^{(1)}} \leq M$. 
Fix $1 \leq i, j \leq d_1$. We have
$$
\begin{aligned}
\mathbb{E}^{x}\Big[ 
\ol{\Xi}_N\big|_{X_i^{(1)}} 
\ol{\Xi}_N\big|_{X_j^{(1)}} \Big]
&= \mathbb{E}^{x}\Big[ 
\ol{\Xi}_N\big|_{X_i^{(1)}} 
\ol{\Xi}_N\big|_{X_j^{(1)}}
-
\ol{\Xi}^{(0)}_N\big|_{X_i^{(1)}} 
\ol{\Xi}^{(0)}_N\big|_{X_j^{(1)}}\Big]
+
\mathbb{E}^{x}\Big[ 
\ol{\Xi}^{(0)}_N\big|_{X_i^{(1)}} 
\ol{\Xi}^{(0)}_N\big|_{X_j^{(1)}}\Big]\\
&=: \mathcal{I}^{(1)}_N+\mathcal{I}^{(2)}_N.
\end{aligned}
$$
It follows from (\ref{harmonicity})
that 
$$
\lim_{N \to \infty}\frac{1}{N}\mathcal{I}_N^{(2)}=
\La \omega_i, \omega_j \Ra_p,
$$
as in the proof of \cite[Lemma 4.2, Step 3]{IKN-1}. 
On the other hand, we obtain
$$
\begin{aligned}
\frac{1}{N}\big|\mathcal{I}^{(1)}_N\big|
&= \frac{1}{N}\Big|\mathbb{E}^x\Big[ 
\big(\Xi_N\big|_{X_i^{(1)}} - \Xi^{(0)}_N\big|_{X_i^{(1)}}\big)
\ol{\Xi}_N\big|_{X_j^{(1)}} \Big]
+\mathbb{E}^x\Big[
\big(\Xi_N\big|_{X_j^{(1)}} - \Xi^{(0)}_N\big|_{X_j^{(1)}}\big)
\ol{\Xi}^{(0)}_N\big|_{X_i^{(1)}} \Big]\Big|\\
&=\frac{1}{N}\Big|\mathbb{E}^x\Big[
\Psi\big(w_N(\cdot)\big)\big|_{X_i^{(1)}}
\ol{\Xi}_N\big|_{X_j^{(1)}} \Big]
+\mathbb{E}^x\Big[
\Psi\big(w_N(\cdot)\big)\big|_{X_j^{(1)}}
\ol{\Xi}^{(0)}_N\big|_{X_i^{(1)}} \Big]\Big|\\
&\leq \frac{M}{N}\Big\{ \mathbb{E}^x \Big[ \big| \ol{\Xi}_N\big|_{X_j^{(1)}} \big|\Big]
+\mathbb{E}^x\Big[ \big| \ol{\Xi}^{(0)}_N\big|_{X_i^{(1)}} \big|\Big]\Big\} \to 0
\end{aligned}
$$
as $N \to \infty$, where we used (\ref{LLN}) for the final line.  \qed

\vspace{2mm}
We write $\alpha^* : \g^{(1)} \to \mathbb{R}$ 
for the {\it Fenchel--Legendre transform} of the quadratic form 
$\alpha(\chi):=\frac{1}{2}\la \Sigma\chi, \chi \ra$ determined by the matrix $\Sigma$, that is, 
$$
\alpha^*(\lambda):=\sup_{\chi \in \g^{(1)}}
\big\{ \la \lambda, \chi \ra
-  \alpha(\chi)\big\}=\sup_{\chi \in \g^{(1)}}
\Big\{ \la \lambda, \chi \ra
-  \frac{1}{2}\la \Sigma\chi, \chi \ra\Big\}
=\frac{1}{2}\la \Sigma^{-1}\lambda, \lambda \ra, \qquad \lambda \in \g^{(1)}. 
$$
Here, $\la \cdot, \cdot \ra$ stand for the standard inner product on 
$\g^{(1)}=\Span\{X_i^{(1)}\}_{i=1}^{d_1} \cong \mathbb{R}^{d_1}$.

Before closing this subsection, we deduce several inequalities which will be used later. 
\begin{lm}\label{estimate above}
Let $\{\ol{\Xi}_n\}_{n=0}^\infty$ and $\{\ol{\Xi}^{(0)}_n\}_{n=0}^\infty$
be as in the previous lemma. 
Then there exist constants $C_1, C_2>0$ independent of $n \in \mathbb{N}$ such that
$$
\|\ol{\Xi}^{(0)}_n(c)\|_{\g^{(1)}} \leq C_1\sqrt{n}, \qquad 
\|\ol{\Xi}_n(c)\|_{\g^{(1)}} \leq C_2\sqrt{n} \qquad \mathbb{P}_x\text{-a.s. }c \in \Omega_x(X).
$$
\end{lm}

\noindent
{\bf Proof.} 
Since $\{\ol{\Xi}^{(0)}_n\}_{n=0}^\infty$ is a (vector-valued) martingale due to the modified harmonicity (cf. \cite[Lemma 3.3]{IKN-1}), we apply the Burkholder--Davis--Gundy inequality to obtain
$$
0 \leq \mathbb{E}^x\big[ \|\ol{\Xi}^{(0)}_n\|_{\g^{(1)}}\big] \leq \mathbb{E}^x\Big[ \sum_{i=1}^n \|\ol{W}_i^{(0)}\|_{\g^{(1)}}^2\Big]^{1/2} \leq C_1\sqrt{n}, \qquad n \in \mathbb{N},
$$
where $\ol{W}_n^{(0)}:=\ol{\Xi}^{(0)}_n - \ol{\Xi}^{(0)}_{n-1}$ for $n \in \mathbb{N}$ and 
$$
C_1^:=\max_{e \in E_0}\Big\|\log\big(\Phi_0(t(\widetilde{e}))\big)\big|_{\g^{(1)}}
-\log\big(\Phi_0(o(\widetilde{e}))\big)\big|_{\g^{(1)}}\Big\|_{\g^{(1)}}<\infty.
$$
We thus conclude the former inequality. As for the latter one, we obtain
$$
\|\ol{\Xi}_n(c)\|_{\g^{(1)}} \leq \|\ol{\Xi}^{(0)}_n(c)\|_{\g^{(1)}} +\max_{x \in V}\|\Psi(x)\|_{\g^{(1)}} 
\leq C_2\sqrt{n}\qquad \mathbb{P}_x\text{-a.s. }c \in \Omega_x(X)
$$
for some $C_2>0$.  \qed
\vspace{2mm}

\subsection{Path spaces} We will deal with several path spaces 
taking its values in $\g^{(1)}$ or $G$ in the next section. 
We denote by $\Map$ the set of all maps $f : [0, 1] \to \g^{(1)}$ satisfying
$f(0)=0$. We endow it with the usual pointwise convergence topology. 
We write $\ContiG$ for the set of all continuous paths $f : [0, 1] \to G$ with 
$f(0)=\bm{1}_G$, unit element in $G$. 
We also denote by $\Conti$ the set of all continuous paths $f : [0, 1] \to \g^{(1)}$ satisfying
$f(0)=0$.
Note that we will consider two different kinds of topologies of $\Conti$
depending on the situations. One is the pointwise convergence topology and
the other is the uniform convergence topology. 
We clearly state which topology is equipped with $\Conti$ in every situation. 
Let $\AC \subset \Conti$ be the subset of all absolutely continuous paths. 
We always equip $\AC$ with the uniform convergence topology.

\section{Proof of Theorem \ref{MDP on G}}

We aim to show the first result (Theorem \ref{MDP on G}) in this section. 
We basically follows the arguments in \cite{Tanaka}, however,
we need to be much careful when those depend on
the speed rate $a_n^2/n$, which is slower than that of LDP-type. 
We may omit the proofs of a few claims if they are done in the same way 
as \cite{Tanaka} and it is not necessary to pay an attention to the speed rate. 
We will use several well-known facts in the general theory of LDPs in this section. 
If readers hope to know details on these facts or their complete proofs, 
we may consult e.g., Dembo--Zeitouni \cite{DZ} and Deuschel--Stroock \cite{DS}.

\subsection{Proof of Proposition \ref{MDP on a path space}}

We define
$$
\ol{Z}^{(n)}_t(c):=\frac{1}{a_n}\ol{W}_{1}(c)+\frac{1}{a_n}\ol{W}_{2}(c)+
\cdots +\frac{1}{a_n}\ol{W}_{[nt]}(c),
\qquad n \in \mathbb{N}, \, 0 \leq t \leq 1, \, c \in \Omega_x(X).
$$
Note that $\{Z^{(n)}\}_{n=1}^\infty$ and $\{\ol{Z}^{(n)}\}_{n=1}^\infty$
are {\it exponentially equivalent} in $\El$ in the sense that 
\begin{equation}\label{exp equiv}
\limsup_{n \to \infty}\frac{n}{a_n^2}\log \mathbb{P}_x\Big(
\|Z^{(n)} - \ol{Z}^{(n)}\|_{\infty} >\delta\Big)=-\infty, \qquad \delta>0,
\end{equation}
where $\El$ is the set of all bounded paths with values in $\g^{(1)}$. 
Indeed, we have
$$
\|Z^{(n)}_t - \ol{Z}^{(n)}_t\|_{\g^{(1)}} \leq \frac{1}{a_n}\|\ol{W}_{[nt]+1}\|_{\g^{(1)}}
\leq \frac{C}{a_n}, \qquad 0 \leq t \leq 1
$$
for some constant $C>0$.  Therefore we obtain
$$
\limsup_{n \to \infty}\frac{n}{a_n^2}\mathbb{P}_x\Big( \|Z^{(n)} - \ol{Z}^{(n)}\|_\infty>\delta\Big)
\leq \limsup_{n \to \infty}\frac{n}{a_n^2}\sum_{k=1}^n 
\mathbb{P}_x\Big( \|\ol{W}_k\|_{\g^{(1)}} > a_n \delta\Big)\leq -\infty
$$
for large $n \in \mathbb{N}$ satisfying $a_n \geq C/\delta$, 
which leads to (\ref{exp equiv}). 
We now set
$$
J:=\{ j=(t_1, t_2, \dots, t_{|j|}) \, | \, 0=t_0<t_1<t_2 <\cdots<t_{|j|} \leq 1\}.
$$
For $j=(t_1, t_2, \dots, t_{|j|}) \in J$, we define
$$
Z^{(n), j}:=\big(Z^{(n)}_{t_1}, Z^{(n)}_{t_2}, \dots, Z^{(n)}_{t_{|j|}}\big) \in (\g^{(1)})^{|j|}, 
\quad
\ol{Z}^{(n), j}:=\big(\ol{Z}^{(n)}_{t_1}, \ol{Z}^{(n)}_{t_2}, \dots, \ol{Z}^{(n)}_{t_{|j|}}\big) \in (\g^{(1)})^{|j|}, 
$$
Then we have the following.
\begin{lm}\label{MDP I_j}
For any $j=(t_1, t_2, \dots, t_{|j|}) \in J$, the sequence of $(\g^{(1)})^{|j|}$-valued
random variables $\{Z^{(n), j}\}_{n=1}^\infty$ satisfies an MDP
with the speed rate $a_n^2/n$ and the good rate function $I_j : (\g^{(1)})^{|j|} \to [0, \infty]$ given by
$$
I_j(\lambda)=\sum_{k=1}^{|j|}(t_k-t_{k-1})\alpha^*\Big(\frac{\lambda_k-\lambda_{k-1}}{t_k-t_{k-1}}\Big),
\qquad \lambda \in (\g^{(1)})^{|j|}.
$$
\end{lm}

\noindent
{\bf Proof.} We will divide the proof into two steps.

\vspace{2mm}
\noindent
{\bf Step~1.} For $j=(t_1, t_2, \dots, t_{|j|}) \in J$, we set
$$
\begin{aligned}
Y^{(n), j}&=\big( \ol{Z}^{(n)}_{t_1}, 
\ol{Z}^{(n)}_{t_2}-\ol{Z}^{(n)}_{t_1}, \dots, \ol{Z}^{(n)}_{t_{|j|}}-\ol{Z}^{(n)}_{t_{|j|-1}}
\big) \in (\g^{(1)})^{|j|}
\end{aligned}
$$
and write $\lambda=(\lambda_1, \lambda_2, \dots, \lambda_{|j|}) \in (\g^{(1)})^{|j|}$. 
Then we have
\begin{align}\label{log-moment}
\mathbb{E}^{x}\Big[ \exp\Big( \frac{a_n^2}{n}\la \lambda, Y^{(n), j}\ra \Big)\Big]
&=\sum_{c \in \Omega_{x, [nt_{|j|}]}(X)}p(c)\exp\Bigg( \frac{a_n}{n}\Big\la \lambda_1, \sum_{i=1}^{[nt_1]}\ol{W}_i(c)\Big\ra + \frac{a_n}{n}\Big\la \lambda_2, \sum_{i=[nt_1]+1}^{[nt_2]}\ol{W}_i(c)\Big\ra
+\cdots\nn\\
&\hspace{1cm}+\frac{a_n}{n}\Big\la \lambda_{|j|}, \sum_{i=[nt_{|j|-1}]+1}^{[nt_{|j|}]}\ol{W}_i(c)\Big\ra\Bigg)\nn\\
&=\sum_{c \in \Omega_{x, [nt_{|j|-1}]}(X)}p(c)\exp\Bigg( \frac{a_n}{n}\Big\la \lambda_1, \sum_{i=1}^{[nt_1]}\ol{W}_i(c)\Big\ra + \frac{a_n}{n}\Big\la \lambda_2, \sum_{i=[nt_1]+1}^{[nt_2]}\ol{W}_i(c)\Big\ra
+\cdots\nn\\
&\hspace{1cm}+\frac{a_n}{n}\Big\la \lambda_{|j|}, \sum_{i=[nt_{|j|-2}]+1}^{[nt_{|j|-1}]}\ol{W}_i(c)
\Big\ra\Bigg)\nn\\
&\hspace{1cm}\times 
\sum_{c' \in \Omega_{t(c), [nt_{|j|}]-[nt_{|j|-1}]}(X)}p(c')
\exp\Big(\frac{a_n}{n}\Big\la \lambda_{|j|}, \ol{\Xi}_{[nt_{|j|}]-[nt_{|j|-1}]}(c')\Big\ra\Big).
\end{align}
Put $N:=[nt_{|j|}]-[nt_{|j|-1}]$ and consider the final expectation 
$$
\mathcal{J}_n(c):=\sum_{c' \in \Omega_{t(c), N}(X)}p(c')
\exp\Big(\frac{a_n}{n}\Big\la \lambda_{|j|}, \ol{\Xi}_N(c')\Big\ra\Big).
$$
By Taylor's formula and the dominated convergence theorem, we have
\begin{align}\label{J_n}
\frac{n}{a_n^2}\log\mathcal{J}_n(c)
&=\sum_{c' \in \Omega_{t(c), N}(X)}p(c')\Bigg\{
 \sum_{i=1}^{d_1}\frac{1}{a_n} 
\ol{\Xi}_N(c')\big|_{X_i^{(1)}}\lambda_{|j|, i}\nn\\
&\hspace{1cm}+\frac{1}{2}\sum_{i, k=1}^{d_1}\frac{1}{n}
\ol{\Xi}_N(c')\big|_{X_i^{(1)}}\ol{\Xi}_N(c')\big|_{X_k^{(1)}}
\lambda_{|j|, i}\lambda_{|j|, k}+O\Big(\frac{a_n}{n^2}\Big)\Bigg\},
\end{align}
where we write $\lambda_{|j|}=(\lambda_{|j|, 1}, \lambda_{|j|, 2}, \dots, \lambda_{|j|, d_1}) \in \g^{(1)}$. 
It follows from (\ref{LLN}), $N \leq n$ and $n/a_n^2 \to 0$ as $n \to \infty$ that
$$
\begin{aligned}
\Big|\sum_{c' \in \Omega_{t(c), N}(X)}p(c')\frac{1}{a_n} 
\ol{\Xi}_N(c')\big|_{X_i^{(1)}}\Big|
&\leq \sum_{c' \in \Omega_{t(c), N}(X)}p(c')\frac{1}{a_n^2} 
\big(\ol{\Xi}_N(c')\big|_{X_i^{(1)}}\big)^2\\
&\leq \frac{n}{a_n^2} \cdot \frac{1}{N}\sum_{c' \in \Omega_{t(c), N}(X)}p(c')
\big(\ol{\Xi}_N(c')\big|_{X_i^{(1)}}\big)^2 \to 0
\end{aligned}
$$
as $n \to \infty$. Combining (\ref{J_n}) with the convergence above and Lemma \ref{convergence of cov}
yields
$$
\begin{aligned}
\lim_{n \to \infty}\frac{n}{a_n^2}\log \mathcal{J}_n(c)
&=\frac{1}{2}(t_{|j|}-t_{|j|-1})\sum_{i, k=1}^{d_1}\La \omega_i, \omega_j \Ra_p 
\lambda_{|j|, i}\lambda_{|j|, k}
=(t_{|j|}-t_{|j|-1}) \cdot \frac{1}{2}\la \Sigma \lambda_{|j|}, \lambda_{|j|} \ra.
\end{aligned}
$$
We repeat the calculation above in (\ref{log-moment}) and we obtain
$$
\begin{aligned}
\Lambda_j(\lambda)&:=\lim_{n \to \infty} \frac{n}{a_n^2}\log \mathbb{E}^x\Big[ \exp\Big( \frac{a_n^2}{n}\la \lambda, Y^{(n), j}\ra \Big)\Big]=\sum_{k=1}^{|j|}(t_k-t_{k-1}) \alpha(\lambda_k)
\end{aligned}
$$ 
for $\lambda=(\lambda_1, \lambda_2, \dots, \lambda_{|j|}) \in (\g^{(1)})^{|j|}$.
Since each $\alpha(\lambda_k) \, (k=1, 2, \dots, |j|)$ is smooth and convex, $\Lambda_j(\lambda)$
is also smooth and convex in $(\g^{(1)})^{|j|}$. Therefore the {\it G\"artner--Ellis theorem} (cf.~\cite[Theorem 2.3.6]{DZ})
implies that
the sequence $\{Y^{(n), j}\}_{n=1}^\infty$ satisfies an MDP
with the speed rate $a_n^2/n$ and the good rate function $\Lambda_j^*(\lambda)$, where
$\Lambda_j^*$ is the Fenchel--Legendre transform of $\Lambda_j$. 
We also write $\chi=(\chi_1, \chi_2, \dots, \chi_{|j|}) \in (\g^{(1)})^{|j|}$. Then we have
$$
\begin{aligned}
\Lambda^*_j(\lambda)&=\sup_{\chi \in (\g^{(1)})^{|j|}}\big\{\la \lambda, \chi \ra - \Lambda_j(\chi)\big\}\\
&=\sum_{k=1}^{|j|}(t_k-t_{k-1})\sup_{\chi \in (\g^{(1)})^{|j|}}
\Big\{ \Big\la \frac{\lambda_k}{t_k-t_{k-1}}, \chi_k\Big\ra - \alpha(\chi_k)\Big\}\\
&=\sum_{k=1}^{|j|}(t_k-t_{k-1})\alpha^*\Big(\frac{\lambda_k}{t_k-t_{k-1}}\Big),
\qquad \lambda \in (\g^{(1)})^{|j|}.
\end{aligned}
$$
Define a continuous map $f : (\g^{(1)})^{|j|} \to (\g^{(1)})^{|j|}$ by 
$$
f(\chi_1, \chi_2, \dots, \chi_{|j|}):=(\chi_1, \chi_1+\chi_2, \dots, \chi_1+\chi_2+\cdots+\chi_{|j|})
$$
for $\chi_k \in \g^{(1)}$ and $k=1, 2, \dots, |j|$. Therefore, by applying the {\it contraction principle}
(cf.~\cite[Theorem 4.2.1]{DZ}), 
the sequence $\{\ol{Z}^{(n), j}=f(Y^{(n), j})\}_{n=1}^\infty$ also satisfies an MDP
with the speed rate $a_n^2/n$ and the good rate function given by
$$
I_j(\lambda):=\sum_{k=1}^{|j|}(t_k-t_{k-1})\alpha^*\Big(\frac{\lambda_k-\lambda_{k-1}}{t_k-t_{k-1}}\Big),
\qquad \lambda \in (\g^{(1)})^{|j|}.
$$

\noindent
{\bf Step~2.} 
By virtue of (\ref{exp equiv}), it is known that $\{Z^{(n), j}\}_{n=1}^\infty$
and $\{\ol{Z}^{(n), j}\}_{n=1}^\infty$ are also exponentially equivalent in $(\g^{(1)})^{|j|}$. 
Indeed, for any $\delta>0$, we have
$$
\begin{aligned}
&\limsup_{n \to \infty}\frac{n}{a_n^2}\log \mathbb{P}_x
\Big( \|Z^{(n), j} - \ol{Z}^{(n), j}\|_{(\g^{(1)})^{|j|}}>\delta\Big)
\leq 
\limsup_{n \to \infty}\frac{n}{a_n^2}\log \mathbb{P}_x
\Big( \|Z^{(n)} - \ol{Z}^{(n)}\|_{\infty}>\delta\Big)=-\infty.
\end{aligned}
$$
Thus the {\it transfer lemma} (cf.~\cite[Theorem 4.2.13]{DZ}) implies that an MDP 
holds for $\{Z^{(n), j}\}_{n=1}^\infty$
with the same good rate function $I_j : (\g^{(1)})^{|j|} \to [0, \infty]$. \qed

\vspace{2mm}
We now see that the set $J$ is viewed as a partially ordered set $(J, \leq)$. 
For $i=(s_1, s_2, \dots, s_{|i|}), \, j=(t_1, t_2, \dots, t_{|j|}) \in J$, we define
$i \leq j$ when $s_\ell=t_{q(\ell)}$ holds  for any $1 \leq \ell \leq |i|$ and for some $1 \leq q(\ell) \leq |j|$. 
Under this partial order, a canonical projection $p_{ij} : M_j \to M_i$ is defined, where 
$M_j:=(\g^{(1)})^{|j|}$ for $j \in J$. 
Let $M$ be the projective limit of the projective system $(M_j, p_{ij})$, that is, 
$M:=\underleftarrow{\lim}\,M_j$. 
Through the canonical projection $p_j : M \to M_j \, (j \in J)$, 
the projective limit $M$ is identified with $\Map$ as follows: 
Since $p_i(f)=p_{ij}\big(p_j(f)\big)$ for $i \leq j$, every $f \in \Map$ corresponds to
$\big(p_j(f)\big)_{j \in J}$. On the contrary, 
every $z=(z_j)_{j \in J} \in M$ is identified with the map $f : [0, 1] \to \g^{(1)}$
given by $f(0)=0$ and $f(t)=z_{\{t\}}$  for $t>0$. 
Moreover, it is known that the topology of $M$ coincides 
with the pointwise convergence topology of $\Map$.
Thanks to the identification above and Lemma \ref{MDP I_j}, 
we may apply the {\it Dawson--G\"artner theorem} (cf. \cite[Theorem 4.6.1]{DZ}).
Therefore, we can show the following by following the same argument as in \cite[Lemma~8.2]{Tanaka}.

\begin{lm} \label{MDP in MAP}
The sequence of $\AC$-valued random variables $\{Z^{(n)}\}_{n=1}^\infty$
satisfies an MDP in $\Map$ equipped with the pointwise convergence topology,
with the speed rate $a_n^2/n$ and the good rate function $I'$ 
defined by {\rm (\ref{rate function path space})}. 
\end{lm}

The rest of key ingredients to show Proposition \ref{MDP on a path space} is 
the following exponential tightness of 
the sequence $\{Z^{(n)}\}_{n=1}^\infty$.  

\begin{lm}\label{exponential tightness}
The sequence of $\AC$-valued random variables 
$\{Z^{(n)}\}_{n=1}^\infty$ is exponentially tight in $\Conti$. 
Namely, for any sufficiently large $R>0$,
there exists a compact set $K_R \subset \Conti$
such that  
$$
\limsup_{n \to \infty} \frac{n}{a_n^2}\log 
\mathbb{P}_x(Z^{(n)} \notin K_R) < -R. 
$$
\end{lm}

\noindent
{\bf Proof.} Since $\|\ol{W}_k\|_{\g^{(1)}}$ has finitely many values 
and is bounded from above for
$k=0, 1, 2, \dots$,
we can choose a sufficiently large $C>0$ such that 
$$
\sup_{k=0, 1, 2, \dots}\max_{i=1, 2, \dots, d_1}\alpha^*(\ol{W}_k|_{X_i^{(1)}}X_i^{(1)}) \leq C.
$$ 
We define a function $\alpha_i^* : \mathbb{R} \to \mathbb{R}$ by
$\alpha^*_i(t):=\alpha^*(t X_i^{(1)})$ for $t \in \mathbb{R}$. 
Then we have $\alpha_i^*(\ol{W}_k |_{X_i^{(1)}}) \leq C$
for $k=0, 1, 2, \dots$ and $i=1, 2, \dots, d_1$. 

For $R>0$, we define a subset $K_R^i \subset \Conti$ by
$$
K_R^i:=\Big\{ f \in \AC \, : \, \int_0^1 \alpha_i^*\big( 
\dot{f}_i(t)\big|_{X_i^{(1)}}\big) \, dt \leq R\Big\}, \qquad i=1, 2, \dots, d_1
$$
and put $K_R:=\bigcup_{i=1}^{d_1}K_R^i$. 
We note that, according to Tanaka \cite[pp.831--832]{Tanaka}, 
we can show that the set $K_R \subset \Conti$ is uniformly bounded and a set of 
equicontinuous functions by using the convexity of each 
$\alpha^*_i : \mathbb{R} \to \mathbb{R} \, (i=1, 2, \dots, d_1)$. 
Then the Ascoli--Arzel\`a theorem leads to the compactness
of $K_R$ in $\Conti$. 

By noting $\dot{Z}^{(n)}_t=(a_n/n)\ol{W}_{[nt]+1}$ for a.e.\,$t \in [0, 1]$, we have
\begin{align}\label{estimate-1}
\mathbb{P}_x(Z^{(n)} \notin K_R)&=\mathbb{P}_x\Big(\bigcup_{i=1}^{d_1}
\{Z^{(n)} \notin K_R^i\}\Big)\nn\\
&\leq d_1 \max_{i=1, 2, \dots, d_1}\mathbb{P}_x(Z^{(n)} \notin K_R^i)\nn\\
&=d_1 \max_{i=1, 2, \dots, d_1}\mathbb{P}_x\Big(\int_0^1 \alpha^*_i
\big(\frac{a_n}{n}\ol{W}_{[nt]+1} \big|_{X_i^{(1)}}\big) \, dt >R\Big)\nn\\
&= d_1 \max_{i=1, 2, \dots, d_1}\mathbb{P}_x\Big(\frac{1}{n}\sum_{k=1}^n \alpha^*_i
\big(\frac{a_n}{n}\ol{W}_{k} \big|_{X_i^{(1)}}\big) >R\Big).
\end{align}
For any $\delta>0$, the Chebyshev inequality yields
\begin{align}\label{estimate-2}
\mathbb{P}_x\Big(\frac{1}{n}\sum_{k=1}^n \alpha^*_i
\big(\frac{a_n}{n}\ol{W}_{k} \big|_{X_i^{(1)}}\big) >R\Big) 
&\leq e^{-\delta R n}
\cdot \mathbb{E}^x
\Big[\exp\Big(\delta \sum_{k=1}^n \alpha^*_i
\big(\frac{a_n}{n}\ol{W}_{k} \big|_{X_i^{(1)}}\big)\Big) \Big]\nn\\
&=e^{-\delta R n}
\cdot\mathbb{E}^x
\Big[\exp\frac{1}{n}\Big(\delta n\sum_{k=1}^n \alpha^*_i
\big(\frac{a_n}{n}\ol{W}_{k} \big|_{X_i^{(1)}}\big)\Big) \Big]\nn\\
&\leq \frac{1}{n}e^{-\delta R \frac{a_n^2}{n} } \sum_{k=1}^n 
\mathbb{E}^x\Big[ \exp\Big(\delta n \alpha^*_i
\big(\frac{a_n}{n}\ol{W}_{k} \big|_{X_i^{(1)}}\big)\Big) \Big],
\end{align}
where we used $e^{-\delta R n} \leq e^{-\delta R \frac{a_n^2}{n}}$ and the Jensen inequality for the convex function $\exp(\cdot)$ for the final line. 
Then it follows from the estimate $\alpha^*_i(\ol{W}_k|_{X_i^{(1)}}) \leq C$ for 
$k=1, 2, \dots, n$ and $i=1, 2, \dots, d_1$ that
\begin{equation}\label{estimate-3}
 \sum_{k=1}^n 
\mathbb{E}^x\Big[ \exp\Big(\delta n \alpha^*_i
\big(\frac{a_n}{n}\ol{W}_{k} \big|_{X_i^{(1)}}\big)\Big) \Big]
\leq n \exp\Big(\delta n \cdot \frac{a_n^2}{n^2} \cdot C\Big)=ne^{C\delta \cdot \frac{a_n^2}{n}}.
\end{equation}
We combine (\ref{estimate-1}) with (\ref{estimate-2}) and (\ref{estimate-3}). Then we have
$$
\begin{aligned}
\frac{n}{a_n^2}\log\mathbb{P}_x(Z^{(n)} \notin K_R)&\leq \frac{n}{a_n^2}\log\Big(
d_1 \cdot \frac{1}{n} \cdot e^{-\delta R \frac{a_n^2}{n}} \cdot ne^{C\delta \frac{a_n^2}{n}}\Big)
=\frac{n}{a_n^2}\log d_1-\delta R + C\delta.
\end{aligned}
$$
Therefore, for a fixed $\delta>0$, we also obtain
$$
\lim_{R \to \infty}\limsup_{n \to \infty}\frac{n}{a_n^2}\log\mathbb{P}_x(Z^{(n)} \notin K_R)=-\infty,
$$
which implies the exponential tightness of the sequence $\{Z^{(n)}\}_{n=1}^\infty$ in $\Conti$. \qed

\vspace{2mm}

We next show Proposition \ref{MDP on a path space}. 

\vspace{2mm}
\noindent
{\bf Proof of Proposition \ref{MDP on a path space}.} 
We note that $\{Z^{(n)}\}_{n=1}^\infty$ is also
regarded as a sequence of random variables with values in $\Conti$
equipped with the pointwise convergence topology and 
the effective domain $\mathcal{D}_{I'}:=\{\varphi \in \Map \, : \, I'(\varphi)<\infty\}$ 
of $I'$ is a subset of $\Conti$.
Since $\{Z^{(n)}\}_{n=1}^\infty$ is exponentially tight in $\Conti$ by Lemma 
\ref{exponential tightness}, 
the {\it inverse contraction principle} 
(cf. \cite[Theorem 4.2.4]{DZ}) implies that $\{Z^{(n)}\}_{n=1}^\infty$
satisfies an MDP in $\Conti$ equipped with the uniform convergence topology
with the good rate function $I'$. Thanks to $\mathcal{D}_{I'} \subset \AC$, 
an MDP holds for the sequence of $\AC$-valued random variables $\{Z^{(n)}\}_{n=1}^\infty$.  \qed

\subsection{Proof of Theorem \ref{MDP on G}}

Define a continuous map 
$F_1 : \AC \to \ContiG$ by 
$F_1(h)(t):=\exp_{F_1(h)(t)}h(t)$ for $h \in \AC$ and $0 \leq t \leq 1$.
Namely, $F_1(h)$ is a solution to the ordinary differential equation $\dot{x}(t)=\dot{c}(t)$ with $x(0)=\bm{1}_G$. 
Note that the continuity of $F_1$ follows from the Gronwall inequality (see \cite[Lemma 43]{Pansu}). 
Let $F_2 : \ContiG \to G$ be a continuous map defined by
$F_2(f):=f(1)$ for $f \in \ContiG$. Put $F:=F_2 \circ F_1 : \AC \to G$, which is also a continuous map. 
Then, by applying the contraction principle, we have the following.

\begin{pr}\label{MDP pre} 
The sequence of $G$-valued random variables $\{F(Z^{(n)})\}_{n=1}^\infty$
satisfies an MDP with the speed rate $a_n^2/n$
and the good rate function $I(g):=\inf \{ I'(h) \, |\, F(h)=g\}$, that is,
$$
\begin{aligned}
-\inf_{g \in A^\circ}I(g) &\leq \liminf_{n \to \infty}\frac{n}{a_n^2}\log\mathbb{P}_x(F(Z^{(n)}) \in A)\\
&\leq  \limsup_{n \to \infty}\frac{n}{a_n^2}\log\mathbb{P}_x(F(Z^{(n)}) \in A)
\leq -\inf_{g \in \ol{A}}I(g)
\end{aligned}
$$
for any Borel set $A \subset G$. 
\end{pr}

By definition, we see that 
$$
F(Z^{(n)})(c)=\exp\Big( \frac{1}{a_n}\ol{W}_1(c)\Big) \cdot 
\exp\Big( \frac{1}{a_n}\ol{W}_2(c)\Big) \cdot \cdots \cdot 
\exp\Big( \frac{1}{a_n}\ol{W}_n(c)\Big), \qquad c \in \Omega_x(X).
$$
We will show that $F(Z^{(n)}(c))$ is close to $\tau_{1/a_n}\big(\ol{\xi}_n(c)\big)$ for sufficiently large $n$
in the sense of the following lemma.

\begin{lm}\label{final lem}
For any $\delta>0$, we have
$$
\limsup_{n \to \infty}\frac{n}{a_n^2}\log \mathbb{P}_x
\Big( d_{\mathrm{Fin}}\big( F(Z^{(n)}), \tau_{1/a_n}(\ol{\xi}_n) \big)>\delta\Big)=-\infty,
$$
where $d_{\mathrm{Fin}}$ is the left invariant {\it Finsler matric} on $G$ defined by
$$
d_{\mathrm{Fin}}(x, y):=
\inf\Big\{ \int_0^1 \|\dot{h}(t)\|_{\g} \, dt \, : \, h \in 
\mathrm{AC}([0, 1]; G), \, h(0)=x, \, h(1)=y\Big\}, \qquad x, y \in G. 
$$
\end{lm}

\noindent
{\bf Proof.} We introduce two $\ContiG$-valued random variables $\phi^{(n)}$ and $\ol{\phi}^{(n)}$ by
$$
\phi^{(n)}(t)(c):=\phi^{(n)}\Big(\frac{k-1}{n}\Big)(c) \exp\Big( 
n\big(t - \frac{k-1}{n}\big) \log \big((\ol{\xi}_{k-1}(c))^{-1} \cdot \ol{\xi}_k(c)\big)\Big)
$$
and 
$$
\ol{\phi}^{(n)}(t)(c):=\ol{\phi}^{(n)}\Big(\frac{k-1}{n}\Big)(c) \exp\Big( 
\frac{n}{a_n}\big(t - \frac{k-1}{n}\big) \ol{W}_k(c)\Big)
$$
for $n \in \mathbb{N}$, $c \in \Omega_x(X)$ and $t \in [(k-1)/n, k/n] \, (k=1, 2, \dots, n)$. Then we see that 
\begin{equation}\label{equalities}
\phi^{(n)}(1)(c)=\ol{\xi}_n(c) \quad \text{and} \quad \ol{\phi}^{(n)}(1)(c)=F(Z^{(n)}(c)), \qquad 
n \in \mathbb{N}, \, c \in \Omega_x(X).
\end{equation}
Moreover, under an appropriate left action, it follows that
$$
(\phi^{(n)})'(t)(c)=n \log \big((\ol{\xi}_{k-1}(c))^{-1} \cdot \ol{\xi}_k(c)\big)
 \quad \text{and} \quad (\ol{\phi}^{(n)})'(t)(c)=\frac{n}{a_n}\ol{W}_k(c) 
$$
for $n \in \mathbb{N}, \, c \in \Omega_x(X)$ and $t \in [(k-1)/n, k/n]\, (k=1, 2, \dots, n)$.
Therefore, for $n \in \mathbb{N}$ and $t \in [(k-1)/n, k/n]\, (k=1, 2, \dots, n)$, we have
$$
\begin{aligned}
&\big\|T_{1/a_n}\big((\phi^{(n)})'(t)(c)\big) - (\ol{\phi}^{(n)})'(t)(c)\big\|_{\g}\\
&\leq \Big\|nT_{1/a_n}\big(\log \big((\ol{\xi}_{k-1}(c))^{-1} \cdot \ol{\xi}_k(c)\big)\big) 
- \frac{n}{a_n}\ol{W}_k(c) \Big\|_{\g}\\
&\leq \frac{n}{a_n^2}\big\|\log \big((\ol{\xi}_{k-1}(c))^{-1} \cdot \ol{\xi}_k(c)\big)
\big|_{\g^{(2)} \oplus \g^{(3)} \oplus \dots \oplus \g^{(r)}}\big\| \leq  C \cdot \frac{n}{a_n^2} \to 0
\end{aligned}
$$
as $n \to \infty$. Then the Gronwall lemma and (\ref{equalities}) imply that, for $c \in \Omega_x(X)$, 
\begin{equation}\label{close}
d_{\mathrm{Fin}}\Big( F(Z^{(n)})(c), \tau_{1/a_n}(\ol{\xi}_n(c)) \Big) \to 0 
\end{equation}
as $n \to \infty$, which leads to the desired convergence. \qed

\vspace{2mm}
\noindent
{\bf Proof of Theorem \ref{MDP on G}.} 
By combining the transfer lemma with Lemma \ref{final lem}, 
the sequence of $G$-valued 
random variables $\{\tau_{1/a_n}(\ol{\xi}_n)\}_{n=1}^\infty$
satisfies an MDP with the speed rate $a_n^2/n$ and 
the good rate function $I : G \to [0, \infty]$, which appears in Proposition \ref{MDP pre}.
Since the map $\varphi : G \to G_\infty$ is a diffeomorphism, 
we use the contraction principle to obtain that 
the sequence of $G_\infty$-valued random variables 
$\{\tau_{1/a_n}\big(\varphi(\ol{\xi}_n)\big)\}_{n=1}^\infty$ satisfies
an MDP with the speed rate $a_n^2/n$ and 
the good rate function $I_\infty : G_\infty \to [0, \infty]$ defined by 
$I_\infty(g):=I\big(\varphi^{-1}(g)\big)$ for $g \in G_\infty$.  \qed


\section{Proof of Theorem \ref{LIL}}

We introduce a compact set of $\AC$ by
\begin{equation}\label{limit points}
\K:=\{ h \in \AC \,  :\, I'(h) \leq 1\},
\end{equation}
where $I'$ is the rate function defined by (\ref{rate function path space}). 
For $h \in \AC$, we define 
$
\mathrm{dist}(h, \K):=\inf \{\|h-h'\|_\infty \, : \, h' \in \K\}.
$ 
The following lemma plays a crucial role to establish Theorem \ref{LIL}. 

\begin{lm}\label{LIL pre}
The sequence $\{Z^{(n)}\}_{n=1}^\infty$ is $\mathbb{P}_x$-a.s.\,relatively compact in
$C_0([0, 1]; \g^{(1)})$. 
\end{lm}

\noindent
{\bf Proof.} We will divide the proof into three steps. 

\vspace{2mm}
\noindent
{\bf Step~1.} We show that, for $\lambda \in \mathbb{N}$,
\begin{equation}\label{toshow1}
\lim_{m \to \infty} \mathrm{dist}(Z^{(\lambda^m)}, \K)=0 \quad \mathbb{P}_x\text{-a.s.}
\end{equation}
For any $\ve>0$, we put $\mathcal{K}_\ve:=\{h \in \AC \, | \, \mathrm{dist}(h, K) \geq \ve\}$. 
Note that, for every $\ve>0$,  the set $\mathcal{K}_\ve$ is compact 
and there exists $\delta=\delta(\ve)>0$ such that 
$
\inf_{h \in \mathcal{K}_\ve}I'(h) > 1+\delta,
$
by the lower semicontinuity of $I'$. Then we have
$$
\begin{aligned}
\sum_{m=1}^\infty \mathbb{P}_x\Big(\mathrm{dist}(Z^{(\lambda^m)}, \mathcal{K})>\ve\Big)
&=\sum_{m=1}^\infty\mathbb{P}_x(Z^{(\lambda^m)} \in \mathcal{K}_\ve)\\
&\leq \sum_{m=1}^\infty\exp\Big( -\inf_{h \in \mathcal{K}_\ve}I'(h) \cdot \log \log \lambda^m \Big)\\
&\leq (\log \lambda)^{-(1+\delta)} \sum_{m=1}^\infty \frac{1}{m^{1+\delta}} <\infty,
\end{aligned}
$$
where we used the moderate deviation estimate with the speed rate $b_n^2/n=\log \log \lambda^m$ 
in Proposition \ref{MDP on a path space} 
for the second line.
Therefore, the first Borel--Cantelli lemma leads to the desired almost sure convergence (\ref{toshow1}). 

\vspace{2mm}
\noindent
{\bf Step~2.} 
For any $n \in \mathbb{N}$, we choose $m \in \mathbb{N}$ and $\lambda \in \mathbb{N}$ such that
$\lambda^m \leq n < \lambda^{m+1}$. We then have
\begin{equation}\label{toshow2}
\mathrm{dist}(Z^{(n)}, \mathcal{K}) \leq 
\mathrm{dist}(Z^{(\lambda^m)}, \mathcal{K}) + \|Z^{(n)} - Z^{(\lambda^m)}\|_\infty.
\end{equation}
It follows from Lemma \ref{estimate above} and 
the time-homogenuity of $\{\ol{\Xi}_n\}_{n=0}^\infty$ that, as $m \to \infty$,
$$
\begin{aligned}
&\|Z_t^{(n)}(c) - Z_t^{(\lambda^m)}(c)\|_{\g^{(1)}}\\
&\leq \frac{1}{\sqrt{\lambda^m \log \log \lambda^m}}
\Big(\frac{\sqrt{\lambda^m \log \log \lambda^m}}{\sqrt{n \log \log n}}-1\Big)
\big\| \ol{\Xi}_{[\lambda^mt]}(c)\big\|_{\g^{(1)}}\\
&\hspace{0.5cm}+\frac{1}{\sqrt{n \log \log n}}\big\|\ol{W}_{[\lambda^mt]+1}(c)+\ol{W}_{[\lambda^mt]+2}(c)
+\cdots + \ol{W}_{[nt]}(c)\big\|_{\g^{(1)}}\\
&\hspace{0.5cm}+\frac{nt-[nt]}{\sqrt{n \log \log n}}\big\|\ol{W}_{[nt]+1}(c)\big\|_{\g^{(1)}}
+\frac{\lambda^mt-[\lambda^mt]}{\sqrt{\lambda^m \log \log \lambda^m}}
\big\|\ol{W}_{[\lambda^mt]+1}(c)\big\|_{\g^{(1)}}\\
&\leq \frac{C_2}{\sqrt{\log\log \lambda^m}}
\Big(\frac{\sqrt{\lambda^m \log \log \lambda^m}}{\sqrt{n \log \log n}}-1\Big)
+\frac{C_2}{\sqrt{\log\log \lambda^m}} + \frac{2C_1}{\sqrt{\lambda^m\log\log \lambda^m}}
\to 0 
\end{aligned}
$$
for $0 \leq t \leq 1$ and $\mathbb{P}_x$-a.s.\,$c \in \Omega_x(X)$, so that we have 
$\|Z^{(n)} - Z^{(\lambda^m)}\|_\infty \to 0$, $\mathbb{P}_x$-a.s. as $m \to \infty$. 
By combining this convergence with (\ref{toshow1}) and (\ref{toshow2}), we obtain
\begin{equation}\label{toshow3}
\lim_{n \to \infty} \mathrm{dist}(Z^{(n)}, \K)=0 \quad \mathbb{P}_x\text{-a.s.}
\end{equation}

\vspace{2mm}
\noindent
{\bf Step~3.} Since $\mathcal{K}$ is compact in $C_0([0, 1]; \g^{(1)})$, 
for any $\ve>0$, we can choose
$\varphi_1, \varphi_2, \dots, \varphi_m \in C_0([0, 1]; \g^{(1)})$ such that 
$\mathcal{K} \subset \bigcup_{i=1}^m B(\varphi_i, \ve/2)$. 
Here, $B(\varphi_i, \ve/2)$ stands for the open ball in $C_0([0, 1]; \g^{(1)})$
centered at $\varphi_i$ of radius $\ve/2$. 
On the other hand, by the previous step, 
there exists $n_0=n_0(\ve) \in \mathbb{N}$ such that $\mathrm{dist}(Z^{(n)}, \mathcal{K})<\ve/2, \, \mathbb{P}_x\text{-a.s.}$ for $n \geq n_0$.  This leads to
$\{Z^{(n)}\}_{n=n_0}^\infty \subset \bigcup_{i=1}^m B(\varphi_i, \ve)$. 
Since $\{Z^{(n)}\}_{n=1}^{n_0} \subset \bigcup_{n=1}^{n_0}B(Z^{(n)}, \ve)$, 
we have
$$
\{Z^{(n)}\}_{n=1}^\infty \subset \Big(\bigcup_{n=1}^{n_0}B(Z^{(n)}, \ve)\Big) \cup 
\Big(\bigcup_{i=1}^m B(\varphi_i, \ve)\Big) \quad \mathbb{P}_x\text{-a.s.}
$$
This implies that $\{Z^{(n)}\}_{n=1}^\infty$ is $\mathbb{P}_x$-a.s.\,totally bounded 
in $C_0([0, 1]; \g^{(1)})$ and therefore
$\{Z^{(n)}\}_{n=1}^\infty$ is $\mathbb{P}_x$-a.s.\,relatively compact 
in $C_0([0, 1]; \g^{(1)})$.  \qed

\vspace{2mm}
The next lemma asserts that an arbitrary element in $\mathcal{K}$ defined by (\ref{limit points})
is a limit point of a suitable subsequence of $\{Z^{(n)}\}_{n=1}^\infty$. 

\begin{lm}\label{LIL-pre2}
The set $\mathcal{K}$ coincides with the set of all 
$\mathbb{P}_x$-a.s.\,limit points of $\{Z^{(n)}\}_{n=1}^\infty$. 
\end{lm}

\noindent
{\bf Proof.} We divide the proof into two steps. 

\vspace{2mm}
\noindent
{\bf Step~1.}
We write $\widetilde{\mathcal{K}}$ for the set of all 
$\mathbb{P}_x$-a.s.\,limit points of $\{Z^{(n)}\}_{n=1}^\infty$. 
For any $h \in \widetilde{\mathcal{K}}$, there exists a subsequence 
$\{Z^{(n_k)}\}_{k=1}^\infty$ of $\{Z^{(n)}\}_{n=1}^\infty$ satisfying 
$\|h - Z^{(n_k)}\|_\infty \to 0$, $\mathbb{P}_x$-a.s.\,as $k \to \infty$. 
We now have
$$
\|h - h'\|_\infty \leq \|h - Z^{(n_k)}\|_\infty + \|Z^{(n_k)} - h'\|_\infty, \qquad h' \in \mathcal{K}.
$$
By taking the infimum running over all $h' \in \mathcal{K}$ and by applying (\ref{toshow3}), we deduce that
$$
\mathrm{dist}(h, \mathcal{K}) \leq \|h - Z^{(n_k)}\|_\infty  + 
\mathrm{dist}(Z^{(n_k)}, \mathcal{K}) \to 0 \quad \text{$\mathbb{P}_x$-a.s.} 
$$
as $k \to \infty$. 
This immediately implies $h \in \mathcal{K}$ and hence we conclude that
$\widetilde{\mathcal{K}} \subset \mathcal{K}$. 

\vspace{2mm}
\noindent
{\bf Step~2.}
On the other hand, we would like to show that $\mathcal{K} \subset \widetilde{\mathcal{K}}$. 
It is sufficient to show that, 
for any $h \in \mathcal{K}$ and $\delta>0$, there exists an integer $\lambda=\lambda(\delta, h)>1$ such that 
\begin{equation}\label{to-show-prob}
\mathbb{P}_x\Big(\limsup_{m \to \infty}
\Big\{ \|Z^{(\lambda^{m})} - h\|_{\infty} \leq \delta\Big\}\Big)=1.
\end{equation}

\noindent
{\bf Step~2.1.}
We introduce linear operators $\T_i \, (i=1, 2, 3)$ on $C([0, 1]; \g^{(1)})$ defined by
$$
\begin{aligned}
\T_1 f(t)&:=f(t) \bm{1}_{[0, 1/\lambda]}(t) + f\Big(\frac{1}{\lambda}\Big)\bm{1}_{(1/\lambda, 1]}(t), &
\T_2 f(t)&:=f\Big(\frac{1}{\lambda}\Big)\bm{1}_{[0, 1/\lambda]}(t) +f(t)\bm{1}_{(1/\lambda, 1]}(t),\\
\T_3 f(t)&:=f\Big(\frac{1}{\lambda^2}\Big) \bm{1}_{[0, 1/\lambda]}(t) + f\Big(\frac{t}{\lambda}\Big)
\bm{1}_{(1/\lambda, 1]}(t) & &
\end{aligned}
$$
for $f \in C([0, 1]; \g^{(1)})$ and $0 \leq t \leq 1$. 
It follows from $f(1/\lambda) \leq \|\T_1 f\|_\infty$ 
and $\|\T_3 f\|_\infty \leq \|\T_1 f\|_\infty$ for $f \in C([0, 1]; \g^{(1)})$ that 
\begin{align}\label{first-est}
\big\|Z^{(\lambda^{m})} - h\big\|_{\infty}
&\leq \|(\T_1+\T_2)Z^{(\lambda^{m})}  - (\T_1+\T_2)h\|_\infty 
+ \|Z^{(\lambda^{m})}_{1/\lambda} - h_{1/\lambda}\|_{\g^{(1)}}\nn\\
&\leq \|\T_1Z^{(\lambda^{m})}  - \T_1h\|_\infty +\|\T_2Z^{(\lambda^{m})}  - \T_2h\|_\infty 
+\|Z^{(\lambda^{m})}_{1/\lambda}\|_{\g^{(1)}} + \|h_{1/\lambda}\|_{\g^{(1)}}\nn\\
&\leq 2\|\T_1Z^{(\lambda^{m})}\|_\infty + 2\|\T_1 h\|_{\infty} + 
\|(\T_2-\T_3)Z^{(\lambda^{m})}  - \T_2h\|_\infty + \|\T_3Z^{(\lambda^{m})}\|_\infty\nn\\
&\leq 3\|\T_1Z^{(\lambda^{m})}\|_\infty + 2\|\T_1 h\|_{\infty} + 
\|(\T_2-\T_3)Z^{(\lambda^{m})}  - \T_2h\|_\infty.
\end{align}
Take an arbitrary $\delta>0$. 
We first see that, thanks to $h_0=0$, there exists a sufficiently large 
$\lambda_1=\lambda_1(\delta, h) \in \mathbb{N}$ such that 
\begin{equation}\label{T_1}
\|\T_1 h\|_\infty < \delta/6
\end{equation}
is true for $\lambda \geq \lambda_1$. 

\vspace{2mm}
\noindent
{\bf Step~2-2.}
We put
$$
F_\delta:=\Big\{ f \in \AC \, : \, \|\T_1 f\|_\infty  \geq \frac{\delta}{9}\Big\},
$$
which is a closed set in $\AC$. 
For $f \in F_\delta$ and $0 \leq t \leq 1/\lambda$, it follows from the Schwarz inequality that 
$$
\begin{aligned}
\frac{\delta}{9} \leq \|f_t\|_{\g^{(1)}}&=\Big\|\int_0^t \dot{f}_s \, ds \Big\|_{\g^{(1)}}
\leq \Big(\int_0^t \|\dot{f}_s\|_{\g^{(1)}}^2 \, ds\Big)^{1/2} \Big(\int_0^{1/\lambda}ds\Big)^{1/2}
=\frac{1}{\sqrt{\lambda}}\Big(\int_0^1 \|\dot{f}_s\|_{\g^{(1)}}^2 \, ds\Big)^{1/2}.
\end{aligned}
$$
Thus we have
\begin{equation}\label{eq1}
\frac{\delta^2}{81}\lambda \leq \int_0^1 \|\dot{f}_s\|_{\g^{(1)}}^2 \, ds. 
\end{equation}
Since the matrix $\Sigma^{-1}$ is positive definite by definition, it holds that 
$\la \Sigma^{-1}\chi, \chi \ra \geq \beta \|\chi\|^2_{\g^{(1)}}$ for $\chi \in \g^{(1)}$,
where $\beta>0$ means the smallest eigenvalue of $\Sigma^{-1}$. 
We then have
\begin{equation}\label{eq2}
\int_0^1 \|\dot{f}_s\|_{\g^{(1)}}^2 \, ds \leq \beta \int_0^1 \la \Sigma^{-1}\dot{f}_s, \dot{f}_s\ra \, ds
\leq 2\beta I'(f).
\end{equation}
By combining (\ref{eq1}) with (\ref{eq2}), we find a sufficiently large 
$\lambda_2=\lambda_2(\delta, h) \in \mathbb{N}$
such that $\lambda \geq \lambda_2$ implies
$$
\inf_{f \in F_{\delta}}I'(f) \geq \frac{\delta^2}{162\beta}\lambda \geq 2.
$$
We next
apply the moderate deviation principle on $\AC$. Namely, 
by the upper bound estimate in Proposition \ref{MDP on a path space}
with the speed rate $b_n^2/n=\log \log \lambda^{m}$, 
we have
$$
\begin{aligned}
\sum_{m=1}^\infty \mathbb{P}_x\Big( \|\T_1 Z^{(\lambda^{m})}\|_\infty  \geq \frac{\delta}{9}\Big)
&=\sum_{m=1}^\infty\mathbb{P}_x(Z^{(\lambda^{m})} \in F_\delta)\\
&\leq \sum_{m=1}^\infty\exp\Big( -\inf_{f \in F_\delta}I'(f) \cdot \log \log \lambda^{m} \Big)\\
&\leq \sum_{m=1}^\infty\exp\big(-2\log \log \lambda^{m}\big) 
\leq \frac{1}{(\log \lambda_2)^2}\sum_{m=1}^\infty\frac{1}{m^2}<\infty
\end{aligned}
$$
for $\lambda \geq \lambda_2$. Thus it follows from the first Borel--Cantelli lemma that 
\begin{equation}\label{T_2}
\mathbb{P}_x\Big( \limsup_{m \to \infty} 
\Big\{ \|\T_1Z^{(\lambda^{m})}\|_\infty \geq \frac{\delta}{9}\Big\}\Big)=0. 
\end{equation}

\vspace{2mm}
\noindent
{\bf Step~2-3.}
Finally, we show that 
\begin{equation}\label{toshowprob}
\mathbb{P}_x\Big( \limsup_{m \to \infty}
\Big\{ \|(\T_2-\T_3)Z^{(\lambda^{m})} - \T_2 h\|_\infty < \frac{\delta}{3}\Big\}\Big)=1.
\end{equation}
We define an open set $G_\delta \subset \AC$ by
$$
G_\delta:=\Big\{ f \in \AC \, : \, 
\|(\T_2-\T_3)f - \T_2 h\|_\infty < \frac{\delta}{3}\Big\}.
$$
Note that $h \in G_\delta$ if we choose $\lambda \geq \lambda_1$. Indeed, we see that 
$$
\|(\T_2-\T_3)h- \T_2 h\|_\infty=\|\T_3 h\|_\infty \leq \|\T_1 h\|_\infty 
< \frac{\delta}{6} < \frac{\delta}{3}, \qquad \lambda > \lambda_1. 
$$
Moreover, it is easy to verify that the sequence $\{(\T_2-\T_3)Z^{(\lambda^{m})}\}_{m=1}^\infty$
is independent since $\{Z^{(\lambda^{m})}\}_{m=1}^\infty$ is independent by definition. 
Then, using the lower bound of the moderate deviation estimate in Proposition 
\ref{MDP on a path space} with the speed rate $b_n^2/n=\log \log \lambda^{m}$ 
and $h \in \mathcal{K}$, we have
$$
\begin{aligned}
&\sum_{m=1}^\infty\mathbb{P}_x\Big(\|(\T_2-\T_3)Z^{(\lambda^{m})} - \T_2 h\|_\infty < \frac{\delta}{3}\Big)\\
&=\sum_{m=1}^\infty\mathbb{P}_x(Z^{(\lambda^{m})} \in G_\delta)\\
&\geq \sum_{m=1}^\infty\exp\Big( -\inf_{f \in G_\delta}I'(f) \cdot \log \log \lambda^{m} \Big)\\
&\geq \sum_{m=1}^\infty\exp\Big( -I'(h) \cdot \log \log \lambda^{m} \Big)
\geq \frac{1}{\log \lambda_1}\sum_{m=1}^\infty\frac{1}{m}=\infty
\end{aligned}
$$
for $\lambda \geq \lambda_1$. Therefore, the second 
Borel--Cantelli lemma leads to (\ref{toshowprob}). 
By combining (\ref{first-est}) with
(\ref{T_1}), (\ref{T_2}) and (\ref{toshowprob}), we obtain (\ref{to-show-prob}) for an arbitrary integer 
$\lambda>\max\{\lambda_1, \lambda_2\}$.  \qed

\vspace{2mm}
For the last, we give a proof of Theorem \ref{LIL}, 
the LIL for the $G_\infty$-valued normalized random walk 
$\{\tau_{1/b_n}(\varphi(\ol{\xi}_n))\}_{n=0}^\infty$. 

\vspace{2mm}
\noindent
{\bf Proof of Theorem \ref{LIL}.} By Lemmas \ref{LIL pre}, \ref{LIL-pre2} and 
the continuity of the map $F : \AC \to G$ defined in Section 3.2, 
we deduce that the sequence $\{F(Z^{(n)})\}_{n=1}^\infty$ is relatively compact 
in $C_{\bm{1}_G}([0, 1]; G)$ and its set of all $\mathbb{P}_x$-a.s.\,limit points coincides with
$$
\widehat{\mathcal{K}}:=F(\mathcal{K})=\{g \in G \, | \, I(g) \leq 1\},
$$ 
where $I : G \to [0, \infty]$ is the good rate function given in Proposition \ref{MDP pre}.  
Let $\widehat{\mathcal{K}}'$ be the set of all $\mathbb{P}_x$-a.s.\,limit points of 
$\{\tau_{1/b_n}(\ol{\xi}_n)\}_{n=0}^\infty$. 
Suppose that $g \in \widehat{\mathcal{K}}$. Then we can choose an appropriate 
subsequence $\{F(Z^{(n_k)})\}_{k=1}^\infty$ satisfying 
$d_{\mathrm{Fin}}(g, F(Z^{(n_k)})) \to 0$, $\mathbb{P}_x$-a.s., as $k \to \infty$. 
By the triangular inequality and (\ref{close}),  we have
$$
d_{\mathrm{Fin}}\big( g, \tau_{1/b_{n_k}}(\ol{\xi}_{n_k})\big)
\leq d_{\mathrm{Fin}}(g, F(Z^{(n_k)})) + d_{\mathrm{Fin}}\big(F(Z^{(n_k)}), 
\tau_{1/b_{n_k}}(\ol{\xi}_{n_k}) \big) \to 0 \quad \mathbb{P}_x\text{-a.s.}
$$
as $k \to \infty$, which implies $g \in \widehat{\mathcal{K}}'$ so that 
$\widehat{\mathcal{K}} \subset \widehat{\mathcal{K}}'$. 
We can also show that $\widehat{\mathcal{K}}' \subset \widehat{\mathcal{K}}$
in the same way as above. 
Therefore we obtain $\widehat{\mathcal{K}}=\widehat{\mathcal{K}}'$. 
Since the map $\varphi : G \to G_\infty$ is a diffeomorphism, 
it follows from Theorem \ref{MDP on G} that 
$$
\widehat{\mathcal{K}}=\{g \in G_\infty \, | \, I(\varphi^{-1}(g)) \leq 1\}
=\{g \in G_\infty \, | \, I_\infty(g) \leq 1\}.
$$
This completes the proof. 
 \qed
 
 \vspace{3mm}
 \noindent
 {\bf Acknowledgement.}
The author would like to thank the anonymous referee for his or her helpful suggestions 
which make the present paper more readable. 
 The author is supported by Grant-in-Aid for JSPS Fellows No.~18J10225 and 
 Grant-in-Aid for Research Activity start-up No.~19K23410.
 



\begin{thebibliography}{9999999}
\bibitem[Ale02]{A2} G.~Alexopoulos:
{\it Random walks on discrete groups of polynomial growth},
Ann.\ Probab.\ {\bf{30}} (2002), pp.~723--801.
\bibitem[BC99]{BC} P. Baldi and L. Caramellino:
{\it Large and moderate deviations for random walks on nilpotent groups},
J. Theor. Probab. {\bf{12}} (1999), pp. 779--809.
\bibitem[BM80]{BM} A.\,A.\,Borovkov and A.\,A.\,Mogulskii:
{\it Probabilities of large deviations in topological space II},
Siberian Math. J. {\bf 21} (1980), pp. 653--664. (Russian original, MR 82c:60049)
\bibitem[CV01]{CV} L. Caramellino and V. Vincenzo:
{\it Law of the iterated logarithm for random walks on nilpotent groups},
Bernoulli {\bf{7}} (2001), pp. 605--628.
\bibitem[CR77]{CR} P. Cr\'epel and B. Roynette: 
{\it Une loi du logarithme it\'er\'e pour le groupe de Heisenberg}, 
Z.~Wahrsch.~Verw.~Gebiete. {\bf 39} (1977), pp. 217--229.
\bibitem[DZ98]{DZ} A. Dembo and O. Zeitouni: 
Large Deviations Techniques and Applications, 2nd Edition, 
Springer, New York, 1998.
\bibitem[DS89]{DS} J.~D.~Deuschel and D.~W.~Stroock: 
Large Deviations, 
Academic Press, New York, 1989.
\bibitem[Goo76]{Goodman} R. W. Goodman: 
Nilpotent Lie Groups, Structure and Applications to Analysis,
LNM {\bf 562}, Springer-Verlag, Berlin Heidelberg, New York, (1976).
\bibitem[Gro81]{Gromov} M. Gromov: 
{\it Groups of polynomial growth and expanding maps}, 
IHES. Publ. Math. {\bf 53} (1981), pp. 53--73.
\bibitem[HL01]{HL} Y. Hu and T. Lee: 
{\it Moderate deviation principles for trajectories of sums of
independent Banach space valued random variables}, 
Trans.~Amer.~Math.~Soc. {\bf 355} (2001), pp. 3047--3064.
\bibitem[IKK17]{IKK} S. Ishiwata, H. Kawabi and M. Kotani:
{\it{Long time asymptotics of non-symmetric random walks on crystal lattices}},
J. Funct. Anal. {\bf 272} (2017), pp.1553--1624.
\bibitem[IKN18-1]{IKN-1} S. Ishiwata, H. Kawabi and R. Namba:
{\it{Central limit theorems for non-symmetric 
random walks on nilpotent covering graphs: Part I}},
Electron. J. Probab. {\bf 25} (2020), 46 pages.  
\bibitem[IKN18-2]{IKN-2} S. Ishiwata, H. Kawabi and R. Namba:
{\it{Central limit theorems for non-symmetric 
random walks on nilpotent covering graphs: Part II}},
Potent. Anal. (2020). Available at {\tt https://doi.org/10.1007/s11118-020-09851-7}. 
\bibitem[Kot02]{Kotani} M. Kotani:
{\it A central limit theorem for magnetic transition operators on a crystal lattice},
J. London Math.\ Soc. {\bf{65}} (2002), pp. 464--482.
\bibitem[Kot04]{Kotani contemp} M. Kotani:
{\it An asymptotic of the large deviation for random walks on a crystal lattice},
Contemp.\ Math.\ {\bf{347}} (2004), pp. 141--152.
\bibitem[KS06]{KS06} M. Kotani and T. Sunada:
{\it Large deviation and the tangent cone at infinity of a crystal lattice},
Math. Z. {\bf{254}} (2006), pp. 837--870.
\bibitem[Mal51]{Malcev} A.~I.~Mal\'cev: {\it On a class of homogeneous spaces}, 
Amer. Math. Soc. Transl. {\bf 39} (1951), pp. 276--307.
\bibitem[Mog76]{Mog} A.~A.~Mogulskii: 
{\it Large deviations for trajectories of multi-dimensional random walks}, 
Theory Probab.~Appl. {\bf 21} (1976), pp. 300--315 (Russian original, MR {\bf 94m}:60063).
\bibitem[Neu96]{Neu} D. Neuenschwander:
{\it Probabilities on the Heisenberg Group. Limit Theorems and Brownian Motion}, 
LNM {\bf 1630}, Berlin, Springer-Verlag, 1996.
\bibitem[Pan83]{Pansu} P. Pansu:
{\it Croissance de boules et des g\'eod\'esiques ferm\'ees dans les nilvari\'et\'es}, 
Ergod. Th. \& Dynam. Sys. {\bf 3} (1983), pp. 415--445.
TEV, Vilnius, 1994.
\bibitem[Sun13]{S} T. Sunada:
Topological Crystallography with a View Towards Discrete  Geometric 
Analysis, Surveys and Tutorials in the Applied Mathematical Sciences {\bf{6}},
Springer Japan, 2013.
\bibitem[Tan11]{Tanaka}
R. Tanaka: 
{\it{Large deviation on a covering graph with group of polynomial growth}}, 
Math. Z. {\bf{267}} (2011), pp. 803--833.
\bibitem[VSC92]{VSC}
N.T. Varopoulos, L. Saloff-Coste and T. Coulhon: 
Analysis and Geometry on Groups,
Cambridge Tracts in Mathematics {\bf{100}}, 
Cambridge Univ. Press, Cambridge, 1992.
\bibitem[Woe00]{Woess} W. Woess: {Random Walks on Infinite Graphs and Groups},
Cambridge Tracts in Mathematics {\bf{138}}. 
Cambridge University Press, Cambridge, 2000. 
\end{thebibliography}
\end{document}